\documentclass[11pt]{article}

\usepackage{amsmath}
\usepackage{amssymb}
\usepackage{latexsym}
\usepackage{amsmath, amsfonts,amssymb, amsthm, euscript,makeidx,color,mathrsfs}

\newtheorem{assumption}{Assumption}
\def\qed{ \ \vrule width.2cm height.2cm depth0cm\smallskip}

\newcommand{\hP}{\hat\dbP}
\newcommand{\res}{{-^{\negmedspace\centerdot\,}}}

\newcommand{\ba}{\begin{array}}
\newcommand{\ea}{\end{array}}
\newcommand{\be}{\begin{equation}}
\newcommand{\ee}{\end{equation}}
\newcommand{\bea}{\begin{eqnarray}}
\newcommand{\eea}{\end{eqnarray}}
\newcommand{\beaa}{\begin{eqnarray*}}
\newcommand{\eeaa}{\end{eqnarray*}}

\def\neg{\negthinspace}

\def\a{\alpha}

\def\e{\varepsilon}

\def\l{\lambda}
\def\m{\mu}

\def\si{\sigma}

\def\o{\omega}

\def\vf{\varphi}

\def\D{\Delta}

\def\O{\Omega}
\def\F{\Phi}

\def\D{\Delta}

\def\F{\Phi}
\def\O{\Omega}
\def\cA{{\cal A}}
\def\cB{{\cal B}}

\def\cE{{\cal E}}
\def\cF{{\cal F}}
\def\cG{{\cal G}}
\def\cH{{\cal H}}

\def\cJ{{\cal J}}

\def\cR{{\cal R}}

\def\cZ{{\cal Z}}
\def\hC{\mathbb{C}}

\def\hE{\mathbb{E}}
\def\hF{\mathbb{F}}

\def\hK{\mathbb{K}}
\def\hL{\mathbb{L}}

\def\hN{\mathbb{N}}

\def\hP{\mathbb{P}}

\def\hR{\mathbb{R}}

\def\hX{\mathbb{X}}

\def\sA{\mathscr{A}}
\def\sB{\mathscr{B}}
\def\sC{\mathscr{C}}

\def\scF{\mathscr{F}}

\def\sK{\mathscr{K}}
\def\scL{\mathscr{L}}

\def\sP{\mathscr{P}}

\def\no{\noindent}

\def\ss{\smallskip}
\def\ms{\medskip}
\def\bs{\bigskip}
\def\q{\quad}
\def\qq{\qquad}

\def\cd{\cdot}

\def\qed{ \hfill \vrule width.25cm height.25cm depth0cm\smallskip}

\newcommand{\basa}{\begin{assumption}}
\newcommand{\easa}{\end{assumption}}

\newcommand{\ol}{\overline}

\newcommand{\bas}{\begin{assum}}
\newcommand{\eas}{\end{assum}}

 \def\cd{\cdot}

\def\as{\hbox{\rm-a.s.{ }}}

\def\co{\mathop{{\rm co}}}

\def\ol{\overline}

\def\cl{\mbox{\rm cl}}
\def\co{\mbox{\rm co}}

\def\dis{\displaystyle}

\def\1{{\bf 1}}

\def\:{\!:\!}

\DeclareMathOperator{\ext}{ext}

\DeclareMathOperator{\dec}{dec}

at 9pt

\newcommand{\of}[1]{\ensuremath{\left( #1 \right)}}
\newcommand{\cb}[1]{\ensuremath{ \left\{ #1 \right\} }}
\newcommand{\sqb}[1]{\ensuremath{ \left[ #1 \right] }}
\newcommand{\norm}[1]{\ensuremath{ \left\Vert #1 \right\Vert }}
\newcommand{\ip}[1]{\ensuremath{ \left\langle #1 \right\rangle }}

\def\prehp(#1,#2){\ensuremath{  #1 \cdot #2 }}

\let\abs=\envert

\begin{document}

\newtheorem{thm}{Theorem}[section]
\newtheorem{lem}[thm]{Lemma}
\newtheorem{cor}[thm]{Corollary}
\newtheorem{prop}[thm]{Proposition}
\newtheorem{rem}[thm]{Remark}
\newtheorem{eg}[thm]{Example}
\newtheorem{defn}[thm]{Definition}
\newtheorem{assum}[thm]{Assumption}

\renewcommand {\theequation}{\arabic{section}.\arabic{equation}}
\def\thesection{\arabic{section}}

\title{\bf  Set-Valued Backward Stochastic Differential Equations}

\author{
 \c{C}a\u{g}{\i}n Ararat\thanks{\noindent Department of Industrial Engineering, Bilkent University, Ankara, 06800, Turkey. E-mail: cararat@bilkent.edu.tr. This author is supported in part by Turkish NSF (T{\"{U}}B\.{I}TAK) 3501-CAREER project \#117F438. The author acknowledges the additional support of the University of Southern California during a research visit for this work in January 2019.}, ~Jin Ma\thanks{ \noindent Department of
Mathematics, University of Southern California, Los Angeles, CA, 90089, USA.
Email: jinma@usc.edu. This author is supported in part by US NSF grant \#1106853. 
},  ~ and ~Wenqian Wu\thanks{\noindent
Department of Mathematics, University of Southern California, Los
Angeles, CA, 90089, USA. E-mail: wenqian@usc.edu. 
 }}

\date{\today}
\maketitle

\begin{abstract}
In this paper, we establish an analytic framework for studying {\it set-valued backward stochastic differential equations (set-valued BSDE)}, motivated largely by the current studies of dynamic set-valued risk measures for multi-asset or network-based financial models. Our framework will  make use of  the notion of {\it Hukuhara difference} between sets, in order to compensate the lack of ``inverse" operation of the traditional Minkowski addition, whence the vector space structure in set-valued analysis. While proving the well-posedness of a class of set-valued BSDEs, we shall also address some fundamental issues regarding generalized Aumann-It\^o integrals, especially when it is connected to the martingale representation theorem. In particular, we propose some necessary extensions of the integral that can be used to represent set-valued martingales with non-singleton initial values. This extension turns out to be essential for the study of set-valued BSDEs.
\end{abstract}

\vfill \bs

\no


\bs

\no{\it 2020 AMS Mathematics subject classification:} 60H05,10; 60G44; 28B20; 47H04.


\eject

\section{Introduction}

Set-valued analysis, both deterministic and stochastic, has found many applications over the years. Most of these applications are in optimization and optimal control theory, but recently more applications have been studied in economics and finance. 
The problem that particularly motivated this work is the so-called {\it set-valued dynamic risk measures}, which we now briefly describe. 

The risk measure of a financial position $\xi$ at a specific time $t$, often denoted by $\rho_t(\xi)$, is defined as a convex functional of 
the (bounded) real-valued random variable $\xi$ satisfying certain axioms such as monotonicity and translativity (cash-additivity) (cf., e.g., \cite{ADEH, bionnadal, gianin}). A {\it dynamic risk measure} is a  family of risk measures $\{\rho_t(\cd)\}_{t\in[0, T]}$, such that for each 
financial position $\xi$, $\{\rho_t(\xi)\}_{t\in[0,T]}$ is an 
adapted stochastic process  satisfying the so-called {\it time-consistency}, in the sense that the following ``tower property" holds (cf. \cite{bionnadal,coquet,gianin}):
\bea
\label{time-consis}
\rho_s(\xi) = \rho_s(-\rho_t(\xi)), \qq \xi\in\hL_{\cF_T}^\infty(\O,\hR),\ 0\leq s\leq t\leq T,
\eea
where $\hL_{\cF_T}^\infty(\O,\hR)$ is the space of $\cF_T$-measurable essentially bounded random variables with values in $\hR$. 
A monumental result in the theory of dynamic risk measures is that, any coherent or even convex risk measure satisfying certain ``dominating" conditions can be represented as the solution of the following {\it 
	Backward Stochastic Differential Equation} (BSDE):
\bea
\label{rmBSDE}
\rho_t(\xi) = -\xi + \int_t^T g(s, \rho_s(\xi),Z_s )ds - \int_t^T Z_s dB_s, \quad t\leq T,
\eea 
where $g$ is determined completely by the properties of  $\{\rho_t\}_{t\in[0,T]}$ (cf. \cite{coquet, MYao, gianin}). 

There has been a tremendous effort to extend the univariate risk measures to the case when the risk appears in the form of a random vector $\xi=(\xi_1,\ldots,\xi_d)\in \hL_{\cF_T}^\infty(\O,\hR^d)$ with $d\in\hN$, typically known as the {\it systemic risk} in the context of default contagion (see, e.g., \cite{conical} for another application in the context of multi-asset markets with transaction costs).
%
%
%
For example, one can consider the contagion of (default) risks in a financial market with large number of institutions  as a {\it network}, in which each institution's future asset value can be viewed as a  ``random shock", to be assessed by its ability to meet its obligations to other members of the network. As a result, it is natural to evaluate these random shocks collectively, which leads to
a multivariate setting of a risk measure, often referred to as ``systemic risk measures" (cf., e.g., \cite{systrisk,fouque,syst}).

One way to characterize a systemic risk measure is to consider it as a multivariate but scalar-valued function. In a static framework, one can 
define an aggregation function $\Lambda\colon \hR^d\to\hR$, so as to essentially reduce the problem to a one-dimensional risk measure. For example, a systemic risk measure can be defined as (cf. \cite{fouque})
\bea
\label{insen}
\rho^{\mathrm{sys}}(\xi)=\rho(\Lambda(\xi))=\inf \{k \in \hR\colon \Lambda(\xi)+k \in \mathcal{A}\},
\eea
where $\xi \in\hL_{\cF_T}^\infty(\O,\hR^d)$ is the wealth vector of the institutions, $\cA$ is a certain {\it acceptance set}, and $\rho$ is a standard
risk measure. 
Such a definition of a systemic risk measure is convenient but have some fundamental deficiencies, especially when one seeks 
a dynamic version. For example, it would be almost impossible to define the tower property (\ref{time-consis}), due to the mis-match of 
the dimensionality. Furthermore, in practice one is often interested in the individual contribution of each institution, and assess the risk
for each institution, thus a more ideal way would be to allocate risks individually, so that the value of a systemic risk measure is defined as a set of vectors. 

It is worth noting that the set-valued risk measure
for a random vector $\xi\in \hR^d$ ($d\geq 2$) can no longer be defined as the ``smallest" capital requirement vector, as it may not exist, for instance, with respect to the componentwise ordering of vectors. One remedy is to 
define it as the set $R_0(\xi)$ (say, at $t=0$) of all the risk compensating portfolio vectors of $\xi$ so that the risk measure $R_0$ is a set-valued functional (see, e.g., \cite{FR2015}). Similarly, one can also define a dynamic set-valued risk measure $\{R_t\}_{t\in[0,T]}$. The 
tower property (\ref{time-consis}) can be defined by 
\bea
\label{tower}
R_s(\xi)=\bigcup_{\eta\in R_t(\xi)}R_s(-\eta)=: R_s[-R_t(\xi)], \qq 0\le s\le t\le T.
\eea
However, the availability of a BSDE-type mechanism to construct or characterize time-consistent dynamic risk measures as in the univariate case is a widely open problem, and is the main purpose of this paper. 

The theory of \emph{set-valued stochastic differential equations (set-valued SDE)} and the related stochastic analysis is not new. Measurability and integration of set-valued functions can
be traced back to as early as 1960s. The commonly used notion of integral is provided by the celebrated work of Aumann \cite{aumann}, where the (Aumann) integral of a set-valued function is defined as the set of all (Lebesgue) integrals of its integrable selections. On the other hand, stochastic integrals of set-valued functions (with respect to Brownian motion or other semimartingales) are relatively new in the literature (see \cite{stochint}). The theory of set-valued SDEs, whose solutions are set-valued stochastic processes (as opposed to \emph{stochastic differential inclusions (SDI)}, whose solutions are vector-valued processes), was established recently (cf., e.g., \cite{KisMic16,MalMic13}). 
While the Backward SDIs have been around for some time (see, e.g., \cite{bsdi,bsdiweak}), to the best of our knowledge, the systematic study of the set-valued BSDEs, especially  in the  general form:
\bea
\label{BSDE-a}
Y_t = \xi + \int^T_t f(s, Y_s, Z_s)ds-\int_t^T Z\circ dB, \qq t\in[0,T],
\eea
is still widely open. (Here, $\int Z\circ dB$ is the {\it generalized} set-valued stochastic integral, see \S3). 

We should point out that the first major difficulty for set-valued analysis, particularly, for studying set-valued BSDEs, is the lack of vector space structure. More precisely, the (Minkowski) addition for sets does not have an ``inverse" (e.g., $A+(-1)A\neq 0$(!)). Consequently, 
 even in the simple case when $f$ is free of $Z$, the equivalence of the BSDE (\ref{BSDE-a}) and its more popular form (cf., e.g., \cite{bsdi, bsdiweak})
\bea
\label{BSDE-b}
Y_t = \hE\Big[\xi + \int^T_t f(s, Y_s)ds\ \Big|\ \cF_t\Big], \qq t\in[0,T],
\eea
is actually not clear at all.

To overcome this difficulty and lay a more generic foundation for the study of BSDEs of type (\ref{BSDE-a}), in this paper we shall explore the notion of the so-called {\it Hukuhara difference} between sets, originated by M. Hukuhara in 1967 \cite{Huku}. We shall first 
establish some fundamental results on stochastic analysis using Hukuhara difference, and then try to prove the the well-posedness of a class of set-valued BSDEs of the form (\ref{BSDE-a}) where $f$ is free of $Z$. It turns out that the seemingly simple additional algebraic structure causes surprisingly subtle technicalities in all aspects of the stochastic analysis, we shall therefore focus on the most basic properties and some key estimates, which will be useful for further development.

Our second goal in this paper is to address some special technical issues in set-valued stochastic analysis involving the generalized Aumman-It\^o integral $\int Z\circ dB$. These issues are subtle, and only occur in the truly set-valued scenarios. When \eqref{BSDE-a} is read as a standard vector-valued BSDE, the indefinite stochastic integral $\int Z_s dB_s = \int Z\circ dB$ appears as a consequence of the classical martingale representation theorem. In the set-valued framework, using the generalized Aumann-It\^o integral, a similar representation theorem was shown in \cite{MgRT} for a set-valued martingale with zero initial value. However, as was pointed out in the recent work \cite{ZhangYano}, if a set-valued stochastic integral is both a martingale and null at zero, then it must be a singleton. Such an observation essentially nullifies any possible role of the martingale representation theorem in the study of set-valued BSDE, unless some modification on the definition of the stochastic integral is adopted. We shall therefore propose a generalization of the Aumman-It\^o
integral so that it contains the information of the non-singleton initial values, and preserves the martingale property. We shall also point out some other fundamental issues regarding the Aumman-It\^o integral in various remarks, but in order not to disturb the main purpose of the paper, we will address these issues in our future publications.

The rest of the paper is organized as follows. In \S2, we give the necessary preliminaries on set-valued analysis, introduce the notion of Hukuhara difference and its properties, and extend the
existing results (mostly in the book \cite{K}) to those that involve Hukuhara difference. In \S3, we revisit set-valued stochastic analysis, again with an eye on these that involve Hukuhara difference. In \S4, we establish some key estimates on set-valued conditional expectations and set-valued Lebesgue integrals. In \S5, we study set-valued martingales and their representations as generalized stochastic integrals. Finally, in \S6, we study the well-posedness of a class of BSDEs of the form (\ref{BSDE-a}) in the case when $f$ is free of $Z$ and compare it to the BSDE of the form (\ref{BSDE-b}).

%
%

\section{Basics of Set-Valued Analysis}
\setcounter{equation}{0}

In this section, we give a brief introduction to set-valued analysis and all the necessary notations associated to it. The interested reader is referred to the books \cite{K, Kis2020} for many of the definitions but we shall present all the results in a self-contained way.

\subsection{Spaces of Sets}

Although most of our discussion applies to more general Hausdorff locally convex topological vector spaces, throughout this paper 
we let $\hX$ be a separable Banach space with norm $\abs{\, \cd\, }$. We shall denote $\sP(\hX)$ to be the set of all nonempty subsets of $\hX$, $\sC(\hX)$ to be the set of all {\it closed} sets in $\sP(\hX)$,  and $\sK(\hX)$  the set of all compact convex sets in $\sP(\hX)$, with respect to the norm topology on $\hX$. We further denote $\sK_w(\hX)$ to be the set of all weakly compact convex sets in $\sP(\hX)$ with respect to the weak topology on $\hX$.

\ms
{\bf  Algebraic Structure on $\sK(\hX)$.} Let $A,B \in \sK(\hX)$ and $\a\in \hR$. We define
\bea
\label{Mink}
A+B := \{a+b\colon a\in A,\ b\in B\}; \q \a A := \{\a a\colon a\in A\}.
\eea
We note that the operations in (\ref{Mink}) are often referred to as the {\it Minkowski addition} and  {\it scalar multiplication}. It can be checked that $\sK(\hX)$ is closed under these operations. It is important to note that the so-called {\it cancellation law} (cf., e.g., \cite{PU,Urbanski}) holds on $\sK(\hX)$, namely, for $A,B,C\in\sK(\hX)$,
\bea
\label{cancellationlaw}
A+C=B+C\q \Longrightarrow\q A=B.
\eea

Clearly, multiplying $A$ by $\a=-1$ gives the ``opposite" of $A$, as $-A:=(-1)A$, which leads to the ``Minkowski difference"
\bea
\label{Minkdiff}
A-B:=A+(-1)B = \{a-b\colon a\in A,\ b\in B\}. 
\eea
But in general, $A+(-1)A \neq \{0\}$, that is, the opposite of $A$ is not the ``inverse" of $A$ under
the Minkowski addition (unless $A$ is a singleton). Consequently, these operations do not establish a vector space structure on $\sK(\hX)$. An early effort to address the inverse operation of Minkowski addition, often still referred to as the {\it Minkowski difference}, is the so-called ``geometric difference" or  ``inf-residuation"  
(see \cite{setoptsurv} and \cite{hamelschrage}), defined by 
\beaa
A\res B := \cb{x\in \hX\mid x+B\subset A},
\eeaa
with $x+B:=\{x\}+B$. Such a difference satisfies $A\res A=\{0\}$, and can be defined for all $A, B\in \sK(\hX)$. However, one 
only has $(A\res B)+B\subset A$; the reverse inclusion usually fails.

In 1967, M.~Hukuhara introduced a definition of set difference that
has since been referred to as {\it Hukuhara difference} (cf. \cite{Huku}) as follows: for $A, B\in \sK(\hX)$,
\bea
\label{Hukudiff}
A \ominus B=C \q\Longleftrightarrow \q A=B+C.
\eea
As we shall see below, this definition has many convenient properties, but the only subtlety is that the Hukuhara difference does not
always exist(!). 
The following result characterizes the existence of Hukuhara difference and  gives an explicit expression of $A\ominus B$, which will be used frequently in our future discussions. 
Recall that, for $A\in\sK(\hX)$ and $a\in A$, $a$ is called an \emph{extreme point} of $A$ if it cannot be written as a strict convex combination of two points in $A$, that is, for every $x_1,x_2\in A$ and $\lambda\in(0,1)$, we have $a\neq \lambda x_1+(1-\lambda)x_2$. 
We denote $\ext(A)$ to be the set of all extreme points of $A$.
\begin{prop}
	\label{Hukuthm}
	Let $A, B\in\sK(\hX)$. The Hukuhara difference $A\ominus B$ exists if and only if for every $a\in\ext(A)$, there exists $x\in \hX$ such that $a\in x+B\subset A$. In this case, $A\ominus B$ is unique, closed, convex, and we have
	\be
	\label{hukuformula}
	A\ominus B = A\res B = \cb{x\in \hX\mid x+B\subset A}.
	\ee
\end{prop}

{\it Proof}: Since this is an infinite-dimensional version of \cite[Proposition~4.2]{Huku} combined with a simple application of the Krein-Milman theorem, we omit the proof.
\qed

The Hukuhara difference facilitates set-valued analysis greatly, without the vector space structure on $\sK(\hX)$. We list some properties that will be used often in this paper.
\begin{prop}
\label{HD0}
Let $A,B,A_1, A_2, B_1, B_2\in \sK(\hX)$, then the following identities hold when all the Hukuhara differences involved exist:

(i) $A\ominus A = \{0\}$, $A \ominus \{0\} = A$;

(ii) $(A_1+B_1) \ominus (A_2+B_2) = (A_1 \ominus A_2) + (B_1 \ominus B_2)$;

(iii) $(A_1+B_1) \ominus B_2 = A_1 + (B_1 \ominus B_2) = (A_1 \ominus B_2) +B_1$;

(iv) $A_1+(B_1\ominus B_2)=(A_1\ominus B_2)+B_1$; and

(v) $A=B+(A\ominus B)$.
\end{prop}

{\it Proof}: (i) $A\ominus A=\{0\}$ is immediate since $A=A+\{0\}$. Suppose $X:=A \ominus \{0\} $. Then by  definition (\ref{Hukudiff}), $A = \{0\} + X = X$. 

(ii) Denote $X:=(A_1+B_1) \ominus (A_2+B_2) $,  $Y:=A_1 \ominus A_2$, and $Z:=B_1 \ominus B_2$. That is,
\be
\label{A}
A_1+B_1 = A_2+B_2+X;\qq
A_1 = A_2+Y; \qq
B_1 = B_2+Z.
\ee
Adding the last two identities above, 
we get
$A_1+B_1 = A_2+Y+B_2+Z = A_2+B_2+Y+Z$.
Comparing this with the first identity in (\ref{A}) and using the cancellation law \eqref{cancellationlaw}, we see that
$X = Y+Z=(A_1 \ominus A_2)+ (B_1\ominus B_2)$,  proving (ii).

(iii) Let $A_2=\{0\}$ in (ii). By the second equality in (i), we obtain the first equality in (iii).  The second equality in (iii) follows by switching the roles of $A_1$ and $B_1$. 

(iv) Denote $X:=B_1\ominus B_2$ and $Y:=A_1\ominus B_2$. That is, $B_1=X+B_2$ and $A_1=Y+B_2$. Then, $A_1+X=Y+B_2+X=Y+B_1$. This is exactly (iv).

(v) This follows immediately by taking $B_1=B_2=B$ in (iv).
\qed

\ms
{\bf Topological Structure on $\sK(\hX)$.}
We note that since $\hX$ is a locally convex topological vector space under both the strong and weak topologies, both $\sK(\hX)$ and $\sK_w(\hX)$ are closed under the Minkowski addition and scalar multiplication. Moreover, the cancellation law \eqref{cancellationlaw}, Proposition~\ref{Hukuthm} and Proposition~\ref{HD0} are valid for both spaces.

For $A,B \in \sK(\hX)$, let us define $\bar{h}(A,B):=\sup_{a\in A}d(a,B)$, where $d(x,B):=\inf_{b\in B}|x-b|$ for $x\in\hX$. Then, the {\it Hausdorff distance} between $A$ and $B$ is given by
\bea
\label{Hdist}
h(A, B):=\bar{h}(A,B)\vee \bar{h}(B,A)=\inf\{\e>0:A\subset V_\e(B), ~ B\subset V_\e(A)\}, 
\eea
where $V_\e(C):=\{ x\in X:d(x,C)\le \e\}$, $C\in \sK(\hX)$, $\e>0$ 
(cf. \cite[Corollary~1.1.3]{K}). Moreover, $(\sK(\hX),h)$ is a Polish space (cf. \cite[Theorem~II.14]{castaing}).
For $A\in \sK(\hX)$, we define
\bea
\label{normdef}
\|A\| := h(A, \{0\})=\sup\{|a|: a\in A\}.
\eea

%

%


We have the following easy results.
\begin{prop}
\label{norm}
(i) The mapping $\|\cd\|:\sK(\hX)\to \hR_+$ satisfies the properties of a norm.

(ii) If $A, B\in \sK(\hX)$ and $A\ominus B$ exists, then $h(A, B) =  \|A\ominus B\|$. 
\end{prop}

{\it Proof.} (i) Clearly, $\|A\| = 0$ implies $A = \{0\}$, and  for any
$\l\in\hR$ we have $ \|\l A\| = h(\l A, \{0\})= \sup\{|\l y|: y\in A\}= |\l|\sup\{|y|: y\in A\}= |\l|\|A\|$. Finally, the ``triangle inequality", in the sense
that $\|A+B\|\le \|A\|+\|B\|$, is trivial by definition of $\|\cd\|$. 

\ss
(ii) Since $A, B \in\sK(\hX)$, applying the translation invariance property of Hausdorff distance (cf. \cite[Proposition 1.3.2]{LGV}), we see that
\bea
\label{Hausnorm}
\|A \ominus B\| = h(A \ominus B, \{0\}) = h((A\ominus B) + B, \{0\} + B) = h(A, B),
\eea
 whenever $A \ominus B$ exists.
\qed

%


%
\begin{rem}
\label{normrem}
{\rm It should be noted that the fact that $\|\cd\|$  satisfies the properties of a norm {\it does not} imply that $( \sK(\hX), \|\cd\|)$ is a normed space, since $\sK(\hX)$ is not a vector space. It is particularly worth noting that, although the Hausdorff metric is symmetric, 
the identity (\ref{Hausnorm}) does not render 
$(A,B)\mapsto \|A\ominus B\|$ a metric on $\sK(\hR^d)$ in the usual sense, since the existence of $A\ominus B$ by no means implies that of $B\ominus A$. In fact, it can be checked that both $A\ominus B$ and $B\ominus A$ exist if and only if $A$ is a translation of $B$ (i.e., $A=x+B$ for some $x\in\hX$). Nevertheless, the relation in Proposition \ref{norm}-(ii) is useful and sufficient for our purposes. 
\qed}
\end{rem}

\subsection{Set-Valued Measurable Functions and Decomposable Sets}\label{svmappings}

We now consider set-valued functions. Let $(E,\cE,\mu)$ be a finite measure space. If $E$ is a topological space, we take $\cE=\sB(E)$, the Borel $\sigma$-algebra on $E$. 
%
%
We shall make use of the following definition of set-valued ``measurable" function. 
\begin{defn}[{\cite[Definition 1.3.1]{Molchanov}}] 
\label{meas1}
A set-valued function $F\colon E \to \sC(\hX)$ is said to be (strongly) measurable if $\{e \in E\colon F(e) \cap B \neq \emptyset\} \in \cE$ for every closed set $B \subset \hX$. 
\end{defn}

The following {\it selection/representation theorems} for set-valued functions are well-known and will be useful in later sections. We shall denote 
$\cl \{A\}$ to be  the closure of a set $A$. 

\begin{thm}\label{KRN}
Let $F\colon E\to\sC(\hX)$ be a set-valued function.\\
(i) \mbox{\rm (Kuratowski and Ryll-Nardzewski, \cite[Theorem 2.2.2]{Kis2020})}  If $F$ is measurable, then $F$ admits a measurable selection, i.e., there exists an $\cE/\sB(\hX)$-measurable function $f\colon E\to \hX$ such that $f(e)\in F(e)$ for each $e\in E$.\\
(ii) \mbox{\rm (Castaing, \cite[Theorem 2.2.3]{Kis2020})}
$F$ is measurable if and only if there exists a sequence $\{f_n\}_{n=1}^\infty$ of measurable selections of $F$ such that $F(e) = \cl\{f_n(e)\colon n\in\hN\}$, $e \in E$. 
\end{thm}


Let us denote $\hL^0(E,\hX)=\hL^0_{\cE}(E,\hX)$ to be the set of all measurable functions $f\colon E\to\hX$ that are distinguished up to $\mu$-almost everywhere (a.e.) equality. For $p\in[1,+\infty)$, let $\hL^p(E,\hX)=\hL^p_{\cE}(E,\hX)$ be the set of all $f\in \hL^0(E,\hX)$ such that $\int_E|f(e)|^p\m(de)<\infty$. Together with the norm $f\mapsto (\int_E|f(e)|^p\m(de))^{\frac1p}$, the set $\hL^p(E,\hX)$ is a Banach space. For $p\in (1,+\infty)$ and $\hX=\hR^d$, $\hL^p(E,\hX)$ is also reflexive. 

We denote $\scL^0(E, \sC(\hX))=\scL^0_{\cE}(E,\sC(\hX))$ to be the space of all measurable set-valued mappings $F\colon E\to \sC(\hX)$ that are distinguished up to $\mu$-a.e. equality. For $F\in \scL^0(E,\sC(\hX))$, we consider the set
\bea 
S(F):= S_{\cE}(F):= \{f \in \hL^0(E,\hX)\colon f(e) \in F(e)\; \mu\text{-a.e.}\ e\in E\}
\eea
of its measurable selections, which is nonempty by Theorem~\ref{KRN}(i). Moreover, by Theorem~\ref{KRN}(ii), two measurable set-valued functions $F$ and $G$ are identical in $\scL^0(E, \sC(\hX))$ if and only if $S(F)=S(G)$. An interesting and crucial question in set-valued analysis is whether a given set of measurable functions in $\hL^0(E,\hR^d)$  can be seen as the set of measurable selections of a measurable set-valued function.  It turns out that this is a highly non-trivial question, for which the following notion is fundamental. 
\begin{defn}
\label{decom}
A set $V\subset \hL^0(E, \hX)$ is said to be 
{\it decomposable} with respect to $\cE$ if it holds ${\bf 1}_D f_1+{\bf 1}_{D^c} f_2\in V$ for every $f_1, f_2\in V$ and $D\in \cE$. 
\end{defn}
Given a set $V\subset \hL^p(E,\hX)$ with $p\in[1,+\infty)$, we define the {\it decomposable hull} of $V$, denoted by $\dec(V)=\dec_{\cE}(V)$, to be the smallest decomposable subset of $\hL^p(E, \hX)$ containing $V$. It can be checked that $\dec(V)$ precisely consists of functions of the form $f=\sum_{i=1}^m {\bf 1}_{D_i}f_i$, where $\{D_1,\ldots,D_m\}$ is a $\cE$-measurable partition of $E$ with $m\in\hN$ and $f_1,\ldots,f_m\in V$. We shall often consider $\overline{\dec}(V)=\overline{\dec}_{\cE}(V)$, the closure of $\dec(V)$ in $\hL^p(E,\hX)$. It is readily seen that $\overline{\dec}(V)$ is the smallest decomposable and closed subset of $\hL^p(E,\hX)$ containing $V$.

For $p\in[1,+\infty)$ and $F\in \scL^0(E, \sC(\hX))$, we define $S^p(F):=S^p_{\cE}(F):= S(F)\cap  \hL^p(E,\hX)$. It is easy to check that $S^p(F)$ is a closed decomposable subset of $\hL^p(E,\hX)$. But it is possible that $S^p(F)=\emptyset$. We  thus consider the set
\bea
\label{Ap}
\sA^p(E,\sC(\hX)):=\sA^p_{\cE}(E,\sC(\hX)):=\{F\in \scL^0(E,\sC(\hX)): S^p(F)\neq \emptyset\},
\eea
and say that $F$ is $p$-integrable if $F\in\sA^p(E,\sC(\hX))$. By \cite[Corollary~2.3.1]{K}, for $F,G\in\sA^p(E,\sC(\hX))$, $F$ and $G$ are identical if and only if $S^p(F)=S^p(G)$. Moreover, we have the following important theorem.
\begin{thm}[{\cite[Theorem~2.3.2]{K}}] 
\label{decomp}
Let $V$ be a nonempty closed subset of $\hL^p(E, \hX)$, $p\ge 1$. Then, there exists $F\in \sA^p(E,\sC(\hX))$ such that  $V=S^p(F)$ if and only if $V$ is decomposable.
\end{thm}

\subsection{Set-Valued Integrals}\label{svintegrals}

We shall now assume that $\hX=\hR^d$, and define the {\it Aumann integral} of a set-valued function $F\colon E\to \sC(\hR^d)$ through its measurable selections. 

As a preparation, for a function $f\in \hL^1(E,\hR^d)$, we define $I(f):=\int_E f(e)\m(de)$ and, for a set $M \subset \hL^1(E, \hR^d)$, we define $I[M]: = \{I(f) : f\in M\}$. Then, one can check (see \cite[Lemma II.3.9]{K}) that $I[M]$ is a convex subset of $\hR^d$ whenever $M$ is decomposable.  Now, for a set-valued function $F\in\sA^1(E, \sC(\hR^d))$, we define
\bea
\label{Aumint}
\int_{E} F(e)\m(de):=\cl (I[S^1(F)]) = \cl\Big\{\int_E f(e)\m(de)\colon f\in S^1(F)\Big\}.
\eea
Clearly, the ``integral" $\int_{E} F(e)\m(de)$ is a nonempty closed convex set, and is called the {\it (closed version of the) Aumann integral} of $F$.

Let $p\in[1,+\infty)$. For a given $F\in\scL^0(E,\sC(\hR^d))$, we say that it is {\it $p$-integrably bounded} if there exists $\ell \in \hL^p(E,\hR_+)$ such that
$\|F(e)\|=h(F(e), \{0\})\le \ell(e)$ a.e. $e\in E$. Let $\scL^p(E,\sC(\hR^d))=\scL^p_{\cE}(E,\sC(\hR^d))$ be the set of all $p$-integrably bounded set-valued functions in $\scL^0(E,\sC(\hR^d))$. It is readily seen that $\scL^p(E,\sC(\hR^d))\subset \sA^p(E,\sC(\hR^d))$. Moreover, by \cite[Theorem~2.4.1-(ii)]{K}, a set-valued function $F\in\sA^p(E,\sC(\hR^d))$ is $p$-integrably bounded if and only if $S^p(F)$ is a bounded subset of $\hL^p(E,\hR^d)$.  In this case, it is even true that $S^p(F)=S^{p^\prime}(F)=S(F)$ for every $p^\prime\in [1,p]$ (cf. \cite[Proposition~2.1.4]{Molchanov}).  In what follows, we shall consider mostly the cases $p=1$ and $p=2$; and say that $F$ is {\it integrably bounded} if 
$F\in \scL^1(E,\sC(\hR^d))$, and {\it square-integrably bounded} if $F\in \scL^2(E,\sC(\hR^d))$. Clearly, $ \scL^2(E,\sC(\hR^d))\subset  \scL^1(E,\sC(\hR^d))$.

We have the following result on integrably bounded set-valued functions. For a subset $A$ of a vector space, $\co(A)$ denotes the convex hull of $A$.

%
\begin{thm}[{\cite[Theorem~2.3.4]{K}}]
 \label{Aumann-co}
 Let $F\in \scL^1(E,\sC(\hR^d))$. Then,
\beaa 
\int_E F(e)\m(de) = \int_E \co(F(e))\m(de).
\eeaa
\end{thm}

In view of Theorem~\ref{Aumann-co}, in the integrably bounded case, it is enough to consider the Aumann integrals of convex-valued functions. On the other hand, if $F\in\scL^p(E,\sC(\hR^d))$, then it is immediate that $F(e)$ is a bounded (hence compact) set for $\mu$-a.e. $e\in E$. In what follows, we mostly restrict our attention to the case $F:E\to\sK(\hR^d)$ and define the spaces $\sA^p(E,\sK(\hR^d))$, $\scL^p(E,\sK(\hR^d))$, and so on in an obvious manner.

Let $F\in \scL^p(E,\sK(\hR^d))$, $p\ge1$. By \cite[Theorem~2.1.18]{Molchanov}, 
we have 
$S^p(F)=S(F)\in \sK_w(\hL^p(E,\hR^d))$. Moreover, since $I$ is a (weakly) continuous linear mapping on $\hL^p(E,\hR^d)$, $I[S^p(F)]=I[S(F)]$ is a nonempty compact convex set  and one can remove the closure in \eqref{Aumint}, that is,
\beaa
\int_E F(e)\m(de)=I[S(F)]\in\sK(\hR^d).
\eeaa  

The following lemma will be helpful in some later calculations. 
\begin{lem}
	\label{meas}
	Let $F_1, F_2\in \scL^p(E, \sK(\hR^d))$, $p\ge1$. Then, $F_1+F_2\in\scL^p(E, \sK(\hR^d))$ and 
	\bea
	\label{addmin1}
	S(F_1+F_2)=S(F_1)+S(F_2).
	\eea
	Furthermore, if $F_1\ominus F_2$ exists, then  $F_1\ominus F_2\in\scL^p(E, \sK(\hR^d))$. In this case, we have 
	\bea
	\label{addmin2}
	S(F_1\ominus F_2)=S(F_1)\ominus S(F_2).
	\eea
\end{lem} 

{\it Proof.} The relation (\ref{addmin1}) is known (see, e.g., \cite[Lemma~2.4.1]{K}). In particular, $S^p(F_1+F_2)\neq \emptyset$ so that $F_1+F_2\in\sA^p(E,\sK(\hR^d))$. Moreover, since $S^p(F_1+F_2)$ is clearly bounded, 
we have $F_1+F_2\in\scL^p(E,\sK(\hR^d))$ and $S^p(F_1+F_2)=S(F_1+F_2)$.

To see the properties of $F_1\ominus F_2$, we assume that it exists. We first claim that $F_1\ominus F_2$ is measurable. Indeed, for $e\in E$ and $x\in\hR^d$, it is easy to check that (cf. \cite[Proposition~4.16]{setoptsurv}) $x\in F_1(e)\ominus F_2(e)$ holds if and only if there exists  a countable dense set $D\subset \hR^d$ (independent of the choice of $e$) such that
\bea
\label{supportfn}
\ip{w,x}\geq \sup_{x_1\in F_1(e)}\ip{w,x_1}-\sup_{x_2\in F_2(e)}\ip{w,x_2}, \qq w\in D.
\eea 
 In other words, we can write
\bea 
\label{supportint}
F_1(e)\ominus F_2(e) = \bigcap_{w\in D}\{x\in\hR^d: \ip{w,x}\geq \sup_{x_1\in F_1(e)}\ip{w,x_1}-\sup_{x_2\in F_2(e)}\ip{w,x_2}\}.
\eea 
Furthermore,  for each $w\in D$, the mappings  $e\mapsto \sup_{x_1\in F_1(e)}\ip{w,x_1}$, $\sup_{x_2\in F_2(e)}\ip{w,x_2}$ 
are measurable real-valued functions by \cite[Example~14.51]{rockwets}, thus every halfspace-valued mapping inside the intersection in \eqref{supportint} is measurable, thus so is the
countable intersection $F_1\ominus F_2$, thanks to \cite[Proposition~14.11-(a)]{rockwets}.

Next, note that $\norm{F_1(e)\ominus F_2(e)}\leq \norm{F_1(e)}+\norm{F_2(e)}$ for every $e\in E$. Since $F_1,F_2$ are $p$-integrably bounded, we see that
$\norm{F_1(\cd)\ominus F_2(\cd)}\in\hL^p(E,\hR)$ and $F_3:=F_1\ominus F_2$ is $p$-integrably bounded. Finally, since $F_2, F_3\in\scL^p(E,\sK(\hR^d))$ and $F_1=F_2+F_3$, \eqref{addmin1} yields $S(F_1)=S(F_2)+S(F_3)$, which then implies that $S(F_1\ominus F_2)=S(F_3)=S(F_1)\ominus S(F_2)$.
\qed

\section{Set-Valued Stochastic Analysis Revisited}
\setcounter{equation}{0}

In this section, we review some basics of set-valued stochastic analysis, and establish some fine results that will be useful for our
discussion but not covered by the existing literature. Throughout the rest of the paper, we shall consider a given complete, filtered
probability space $(\O, \cF, \hP, \hF=\{\cF_t\}_{t\in[0,T]})$, on which is defined a standard $m$-dimensional Brownian motion $B=\{B_t\}_{t\in[0,T]}$, where $T>0$ is
a given time horizon. We shall denote $\hL^p_{\hF}([0,T]\times \Omega,\hR^d)$ to be the space of all $\hF$-progressively measurable $d$-dimensional processes $\{\phi_t\}_{t\in[0,T]}$ with $\hE[\int_0^T |\phi_t|^p dt]<+\infty$. The space $\hL_{\hF}^p([0,T]\times\O,\hR^{d\times m})$ of matrix-valued processes can be defined similarly.

\subsection{Set-Valued Conditional Expectations}\label{svrv}

A set-valued random variable $X:\O\to \sC(\hR^d)$ is an $\cF$-measurable set-valued function. If $X \in \sA^1(\O, \sC(\hR^d))$, then we define its expectation, denoted by $\hE[X]$ as usual, by its Aumann integral $\int_\O X(\o)\hP(d\o)$. Given $p\ge1$, if $X \in \sA^p(\O, \sC(\hR^d))$, then
$S^p(X)$ is a closed decomposable subset of $\hL^p(\O, \hR^d)$ and $S^p(\ol{\co}(X))=\ol{\co}(S^p(X))$ (see \cite[Lemma 2.3.3]{K}).
Further, $X$ is $p$-integrably bounded if and only if  $S^p(X)$ is a bounded set in $\hL^p(\O,\hR^d)$, that is, 
$\hE[\|X\|^p]=\int_\O\sup\{|x|^p: x\in X(\o)\}\hP(d\o)=\int_\O h^p(X(\o), \{0\})\hP(d\o)<\infty$. In particular, if $X\in \scL^p(\O,\sK(\hR^d))$, then $S^p(X)=S(X)$ is a weakly compact convex subset of $\hL^p(\O,\hR^d)$.

	Let $\cG$ be a sub-$\sigma$-field of $\cF$. We denote
	$\hL^p_{\cG}(\O,\hR^d)$, $\sA_{\cG}^p(\O,\sC(\hR^d))$, $\scL_{\cG}^p(\O,\sK(\hR^d))$, $S^p_{\cG}(X)$ to be the same as those in Sections \ref{svmappings} and \ref{svintegrals}, on the probability space $(\O, \cG, \hP)$.
%
%
%
%
Further, for $X\in \sA_\cF^{1}(\O, \sC(\hR^d))$, the {\it conditional expectation} of $X$ given $\cG$ is defined as the (almost surely) unique set-valued random variable $\hE[X|\cG]\in \sA_\cG^{1}(\O, \sC(\hR^d))$ that satisfies
\bea
\label{set-valued condexp}
S_{\cG}^1(\hE[X|\cG]) = \cl\{\hE[f|\cG]: f\in S^1(X)\},
\eea 
where the closure is evaluated in $\hL_{\cG}^1(\O,\hR^d)$. The existence of $\hE[X|\cG]$ follows by Theorem~\ref{decomp} since the set on the right in \eqref{set-valued condexp} is decomposable. Moreover, for $p\ge 1$, if $X\in\scL_{\cF}^p(\O,\sK(\hR^d))$, then it can be shown that the closure in \eqref{set-valued condexp} is not needed and $\hE[X|\cG]\in\scL_{\cG}^p(\O,\sK(\hR^d))$. In this case, $\hE[X|\cG]$ satisfies the usual identity
\bea
\label{condexp}
\int_D \hE[X|\cG](\o)\hP(d\o)=\int_D X(\o)\hP(d\o), \qq  D\in \cG.
\eea
Moreover, it can be easily checked that $\hE[\cd|\cG]$ satisfies all the natural properties of a conditional expectation, 
except that the ``linearity" should be interpreted in terms of the Minkowski addition and multiplication by scalars. Furthermore, we note that the conditional expectation of a set $V\subset \hL_{\cF}^1(\O,\hR^d)$ of random variables can also be defined in a generalized sense even if it is not the set of selections of a set-valued random variable. To be more precise, if $V\subset\hL_{\cF}^1(\O, \hR^d)$ is a nonempty closed decomposable set, then there exists a unique $\hE[V|\cG]\in \sA_{\cG}^1(\O,\sC(\hR^d))$ (by a slight abuse of notation) such that
\bea
\label{condexpg}
S_{\cG}^1(\hE[V|\cG]) = \cl\{\hE[f|\cG]: f\in V\}.
\eea 

%
%
%
%
%

The following is a seemingly obvious fact regarding set-valued conditional expectations.
\begin{cor}
	\label{condi}
	Let $X_1, X_2\in \scL^p(\O,\sK(\hR^d))$ with $p\in[1,+\infty)$. Let $\cG\subset \cF$ be a sub-$\si$-field. Suppose that $X_1\ominus X_2$ exists. Then, $ \hE[X_1\ominus X_2|\cG]$  exists in $\scL_{\cG}^p(\O,\sK(\hR^d))$ and it holds that
	\bea
	\label{condhuku}
	\hE[X_1\ominus X_2|\cG]=\hE[X_1|\cG]\ominus \hE[X_2|\cG].
	\eea
\end{cor}

{\it Proof.} By Lemma~\ref{meas}, $X_1\ominus X_2\in\scL^p(\O,\sK(\hR^d))$ so that $\hE[X_1\ominus X_2|\cG]$ exists in $\scL_{\cG}^p(\O,\sK(\hR^d))$. By the definition of conditional expectation and repeated applications of Lemma~\ref{meas}, we have
\beaa
S_{\cG}(\hE[X_1\ominus X_2|\cG]+\hE[X_2|\cG])
&=&S_{\cG}(\hE[X_1\ominus X_2|\cG])+S_{\cG}(\hE[X_2|\cG])\\
&=&\{\hE[f_1|\cG]:f_1\in S(X_1\ominus X_2)\}+\{\hE[f_2|\cG]:f_2\in S(X_2)\}\\
&=&\{\hE[f|\cG]:f\in S(X_1\ominus X_2)+S(X_2)\}\\
&=&\{\hE[f|\cG]:f\in S((X_1\ominus X_2)+X_2)\}=S_{\cG}(\hE[X_1|\cG]).
\eeaa
This is equivalent to having $\hE[X_1|\cG]=\hE[X_1\ominus X_2|\cG]+\hE[X_2|\cG]$, whence \eqref{condhuku}. 
\qed

\subsection{Set-Valued Stochastic Processes}\label{svsp}

 A set-valued stochastic process $\Phi=\{\Phi_t\}_{t\in[0,T]}$ is a family of set-valued random variables taking values in $\sC(\hR^d)$. We call $\Phi$ measurable if it is $\sB([0,T])\otimes\cF$-measurable as a single set-valued function on $[0,T]\times\O$.  The notions such as ``adaptedness" or ``progressive measurability" can be defined accordingly in the obvious ways. We denote $\scL^0_\hF([0,T]\times\O,\sC(\hR^d))$ to be the space of all set-valued, $\hF$-progressively measurable processes taking values in $\sC(\hR^d)$.
For $\Phi\in \scL^0_\hF([0,T]\times\O,\sC(\hR^d))$, we denote $S_{\hF}(\Phi)$ to be the set of all $\hF$-progressively measurable selectors of $\Phi$, which is nonempty by Theorem \ref{KRN}. For $p\in[1,+\infty)$, we define $S^p_{\hF}(\Phi):=S_{\hF}(\Phi)\cap\hL^p_{\hF}([0,T]\times\O,\hR^d)$ and denote $\scL^p_\hF([0,T] \times \Omega, \sC(\hR^d))$ to be the set of all $\hF$-progressively measurable,
$\sC(\hR^d)$-valued processes $\Phi$ with $\hE[\int_0^T\norm{\Phi_t}^p dt]<+\infty$ (i.e., $p$-integrably bounded). The notations $\scL^p_\hF([0,T] \times \Omega, \sK(\hR^d)), \scL_{\hF}^p([0,T]\times\O,\sK(\hR^{d\times m}))$ for set-valued processes with compact convex values are defined similarly for $p=0$ and $p\geq 1$. It is worth pointing out that the space $\scL^2_\hF([0,T] \times \Omega, \sK(\hR^d))$ is not a Hilbert space, but only a complete metric space, with the metric $d_H(\Phi, \Psi):=(\hE[\int_0^T
h^2(\Phi_t, \Psi_t)dt])^{1/2}$. 
%

\subsection{Set-Valued Stochastic Integrals}\label{sec:stochint}

	In this section, we assume that $\hF=\hF^B$, the natural filtration generated by $B$, augmented by all the $\hP$-null sets of $\cF$ so that it satisfies the {\it usual hypotheses}.

Let us consider the two linear mappings $J : \hL^2_\hF([0,T]\times \Omega,\hR^d)\to \hL^2_{\cF_T}(\O,\hR^d)$, and $\cJ : \hL^2_\hF([0,T]\times \O,\hR^{d\times m})\to\hL^2_{\cF_T}(\O, \hR^d)$ defined by
\bea
\label{JsJ}
 J(\phi): = \int^T_0\phi_tdt, \qq \cJ(\psi) :=\int^T_0\psi_tdB_t,
\eea
for $\phi \in \hL^2_\hF([0,T]\times \Omega,\hR^d)$, $\psi \in \hL^2_\hF([0,T]\times \O,\hR^{d\times m})$, respectively.
For $K\subset\hL^2_\hF([0,T]\times \Omega,\hR^d)$ (resp. $K^\prime \subset\hL^2_\hF([0,T]\times \O,\hR^{d\times m})$), the set $J[K]$ 
(resp. $\cJ[K^\prime]$) is defined in an obvious way. 

Let $\Phi\in \scL^0_\hF([0,T]\times\O,\sC(\hR^d))$ and $\Psi \in \scL^{0}_\hF([0,T]\times\O,\sC(\hR^{d\times m}))$ such that $S^2_\hF(\Phi)\neq\emptyset$, $S^2_\hF(\Psi)\neq\emptyset$.
Then, one can show that there exist unique set-valued random variables $\int^T_0\Phi_tdt\in \sA^{2}_{\cF_T}(\O,\sC(\hR^d))$ and 
$\int^T_0\Psi_tdB_t\in \sA^{2}_{\cF_T}(\O,\sC(\hR^d))$ such that 
\bea
\label{itoint}
 S^2_{\cF_T}\of{\int^T_0\Phi_tdt}=\overline{\dec}_{\cF_T}(J[S^2_\hF(\Phi)]), \qq
S^2_{\cF_T}\of{\int^T_0\Psi_tdB_t}=\overline{\dec}_{\cF_T}(\cJ[S^2_\hF(\Psi)]).
\eea
We call $\int_0^T \Phi_t dt$ and $\int_0^T\Psi_t dB_t$  {\it set-valued  stochastic integrals}. As usual, for $t\in[0,T]$, we define the {\it indefinite} stochastic integrals as $\int_0^t \Phi_s ds:=\int_0^T\1_{(0,t]}(s)\Phi_sds$ and $\int_0^t \Psi_s dB_s:=\int_0^T\1_{(0,t]}(s)\Psi_sdB_s$. Equivalently, one can define them via the relations $S^2_{\cF_t}(\int_0^t \Phi_s ds)=\overline{\dec}_{\cF_t}(J_{0,t}[S^2_{\hF}(\Phi)])$, $S^2_{\cF_t}(\int_0^t \Psi_s dB_s)=\overline{\dec}_{\cF_t}(\cJ_{0,t}[S^2_{\hF}(\Psi)])$, where $J_{0,t}(\phi):=\int_0^t \phi_sds$, $\cJ_{0,t}(\psi):=\int_0^t \psi_sdB_s$. The integrals $\int_t^T \Phi_s ds$ and $\int_t^T \Psi_s dB_s$, and the mappings $J_{t,T}$, $\cJ_{t,T}$ can be defined similarly for $t\in[0,T]$.

\begin{rem}
\label{propint}
{\rm The set-valued  It\^o stochastic integrals have
many interesting properties, we refer the interested reader to the books \cite{K, Kis2020} for the exhaustive explorations. Here we mention a few that will 
be useful for our discussion. 

 (i) The definition  \eqref{itoint} implies that both $\int_0^T\Phi_t dt$ and $\int_0^T\Psi_t dB_t$ are $\cF_T$-measurable set-valued random 
 variables. However,  neither of the sets $J[S^2_\hF(\Phi)], \cJ[S^2_\hF(\Psi)]\subset \hL^2_{\cF_T}(\O, \hR^d)$ is necessarily decomposable 
 (see \cite[p.105]{K} for counterexamples). Thus, by virtue of Theorem  \ref{decomp}, they cannot be seen as the selectors of any 
 $\cF_T$-measurable set-valued random variables. 
 
 (ii) One can actually show that $\{\hE[x]\colon x\in \cJ[S^2_\hF(\Psi)]\}=\{0\}$, and $\cJ[S^2_\hF(\Psi)]$ is decomposable if and only if it is a singleton(!).
 
 (iii) By \cite[Theorem 3.1.1]{K}, it is shown that $\overline{\dec}_{\cF_T}(\cJ[S^2_\hF(\Psi)])=\hL^2_{\cF_T}(\O,\hR^d)$ if and only int $ \overline{\dec}(\cJ[S^2_\hF(\Psi)])\neq \emptyset$. 
 
 (iv) If $\Phi$ and $\Psi$ are convex-valued, then so are $\int_0^T \Phi_t dt$ and $\int_0^T\Psi_t dB_t$. If $\Phi\in \scL^2_\hF([0,T]\times\O,\sK(\hR^d))$, then it is known that $\int_0^T \Phi_tdt\in\scL^2_{\cF_T}(\O,\sK(\hR^d))$, that is, the stochastic time integral of a square-integrably bounded process is a square-integrably bounded set-valued random variable (see \cite[Theorem 3.2]{Kproperties}). However, the It$\hat{\text{o}}$ integral $\int_0^T \Psi_tdB_t$ fails to be square-integrably bounded in general even if $\Psi\in \scL^{2}_\hF([0,T]\times\O,\sK(\hR^{d\times m}))$ (see \cite{Munbounded}).
 
(v) The set-valued stochastic integrals $\int_0^t \Phi_sds$, $\int_0^t\Psi_sdB_s$ are defined, almost surely, for each $t\in[0,T]$, and they are ($\hF$-){\it adapted}, in the usual sense. Furthermore, 
when $\Phi\in \scL^{2}_\hF([0,T]\times\O,\sK(\hR^{d}))$, by \cite[Theorem 2.4]{KisMic17}, the process $\{\int_0^t \Phi_sds\}_{t\in[0,T]}$ has a continuous (with respect to $h$), whence progressively measurable, modification. We can define the indefinite integral $\int_0^\cd \Psi_sds$ by this progressively measurable set-valued process. However, the continuity of the It\^o integral $\{\int_0^t \Psi_sdB_s\}_{t\in[0,T]}$ is much more involved, and so is the progressive measurability  issue (see \cite[Section 5.5]{Kis2020} for a special case).
 \qed}
 \end{rem}

The following lemma shows that the
additivity holds for both integrals, which also allows to calculate the integrals of the Hukuhara difference of two processes. 
\begin{lem}
	\label{stochint-add}
	Suppose that $\hP$ is a nonatomic probability measure. Let $\Phi^1,\Phi^2\in\scL^2_{\hF}([0,T]\times\O,\sK(\hR^d))$ and $\Psi^1,\Psi^2\in \scL^2_{\hF}([0,T]\times\O,\sK(\hR^{d\times m}))$. Then, for every $t\in[0,T]$,
	\bea 
	\label{addint}
	\int_0^t (\Phi^1_s+\Phi^2_s)ds = 	\int_0^t \Phi^1_s ds +	\int_0^t \Phi^2_s ds ,\quad \int_0^t (\Psi^1_s+\Psi^2_s)dB_s = 	\int_0^t \Psi^1_s dB_s+	\int_0^t \Phi^2_s dB_s
	\eea 
	hold almost surely.
	If $\Phi^1\ominus \Phi^2$ and $\Psi^1\ominus \Psi^2$ exist ($dt\times d\hP$-a.e.), then we have $\Phi^1\ominus \Phi^2\in\scL^2_{\hF}([0,T]\times\O,\sK(\hR^d))$, $\Phi^1\ominus \Phi^2\in\scL^2_{\hF}([0,T]\times\O,\sK(\hR^{d\times m}))$ and, for every $t\in[0,T]$,
	\bea 
	\label{Hukuharaint}
	\int_0^t (\Phi^1_s\ominus \Phi^2_s)ds = 	\int_0^t \Phi^1_s ds \ominus \int_0^t \Phi^2_s ds ,\quad \int_0^t (\Psi^1_s\ominus \Psi^2_s)dB_s = 	\int_0^t \Psi^1_s dB_s\ominus	\int_0^t \Phi^2_s dB_s
	\eea 
	hold almost surely.
	\end{lem}
 
 {\it Proof}: The relations in \eqref{addint} are given in \cite[Theorem 3.1-3.2]{Kproperties}. Suppose that $\Phi^1\ominus \Phi^2$ exists. It is clear that $\Phi^1\ominus\Phi^2$ takes values in $\sK(\hR^d)$. Since $\norm{\Phi_t^1\ominus \Phi_t^2}\leq \norm{\Phi_t^2}+\norm{\Phi_t^2}$,
 \beaa
 \hE\sqb{\int_0^T \norm{\Phi_t^1\ominus\Phi^2_t}^2 dt }\leq 2\hE\sqb{\int_0^T\norm{\Phi_t^1}^2dt}+2\hE\sqb{\int_0^T\norm{\Phi_t^2}^2dt}<+\infty.
 \eeaa
 This and Lemma \ref{meas} imply that $\Phi^1\ominus \Phi^2\in\scL^2_{\hF}([0,T]\times\O,\sK(\hR^d))$. We have $\Phi^1=\Phi^2+(\Phi^1\ominus\Phi^2)$. Let $t\in[0,T]$. By the first relation in \eqref{addint}, we obtain $\int_0^t \Phi_s^1 ds = \int_0^t \Phi_s^2 ds + \int_0^t (\Phi_s^1\ominus \Phi_s^2)ds$. By the definition of Hukuhara difference, the first relation in \eqref{Hukuharaint} follows. The proofs of the claims related to $\Psi^1\ominus \Psi^2$ are similar, hence omitted.
 \qed 

\begin{cor}
	\label{dtint-time}
	Suppose that $\hP$ is a nonatomic probability measure. Let $\Phi\in\scL^2_{\hF}([0,T]\times\O,\sK(\hR^d))$, $\Psi\in\scL^2_{\hF}([0,T]\times\O,\sK(\hR^{d\times m}))$. For each $t\in [0,T]$,
	\[
	\int_0^T \Phi_s ds = \int_0^t \Phi_s ds + \int_t^T \Phi_s ds,\quad \int_0^T \Psi_s dB_s = \int_0^t \Psi_s dB_s + \int_t^T \Psi_s dB_s
	\]
	and
	\[
	\int_t^T \Phi_sds = \int_0^T \Phi_s ds \ominus \int_0^t \Phi_s ds,\quad \int_t^T \Psi_sdB_s = \int_0^T \Psi_s dB_s \ominus \int_0^t \Psi_s dB_s.
	\]
 hold almost surely.
	\end{cor}

 {\it Proof}: This is immediate from Lemma \ref{stochint-add} and the definitions of the integrals since $\1_{(0,T]}(s)\xi_s = \1_{(0,t]}(s)\xi_s +\1_{(t,T]}(s)\xi_s$, $\xi\in\{\Phi, \Psi\}$, 
 for all $s\in[0,T]$.
 \qed

The notion of stochastic integral can be extended to the case where the integrand is only a set of processes, instead of a set-valued process. We briefly describe the idea (cf. \cite{MgRT}). Let $\cZ\in \sP(\hL^2_\hF([0,T]\times\O,\hR^{d\times m}))$ be a nonempty set and consider 
the sets $\cJ_t[\cZ]=\{\int_0^t z_sdB_s: z\in \cZ\}$, $t\in[0,T]$. Due to lack of decomposability, $\cJ_t[\cZ]$ is not equal to the set of square-integrable selections 
of a set-valued random variable, in general. But similar to 
the stochastic integral discussed above, one can show that,  for each $t\in[0,T]$, there exists a unique $\int_0^t \cZ\circ dB\in \sA^2_{\cF_t}(\O,\sC(\hR^d))$ such that
\bea\label{genint}
S^2_{\cF_t}\of{\int_0^t\cZ\circ dB}= \overline{\dec}_{\cF_t}(\cJ_t[\cZ]).
\eea 
We call $\int_0^t\cZ\circ dB$ the {\it generalized (indefinite) Aumann-It\^o stochastic integral} (cf. \cite{MgRT}). If $\cZ$ is convex, then $\int_0^t \cZ\circ dB$ is convex-valued (see \cite[Theorem 2.2]{MgRT}). 

We have the following analogue of Lemma \ref{stochint-add}.

\begin{lem}
	\label{genstochint-add}
	Assume that $\hP$ is nonatomic, and let $\cZ^1,\cZ^2\in\sK_w(\hL^2_{\hF}([0,T]\times\O,\hR^{d\times m}))$. Then, 
	the following statements are true:
	
	(i) $\cZ^1+ \cZ^2\in \sK_w(\hL^2_{\hF}([0,T]\times\O,\hR^{d\times m}))$ and for every $t\in[0,T]$, it holds that 
	\bea
	\label{add1} 
	\int_0^t (\cZ^1+\cZ^2)\circ dB =\int_0^t \cZ^1\circ dB + \int_0^t \cZ^2\circ dB, \qq \hP\as
	\eea

 (ii) If $\cZ^1\ominus \cZ^2$ exists, then $\cZ^1\ominus \cZ^2\in \sK_w(\hL^2_{\hF}([0,T]\times\O,\hR^{d\times m}))$ and for every $t\in[0,T]$,
	\bea
	\label{sub1}
	\int_0^t (\cZ^1\ominus \cZ^2)\circ dB =\int_0^t \cZ^1\circ dB \ominus \int_0^t \cZ^2\circ dB, \qq \hP\as
	\eea

(iii) If $\cZ^1\ominus \cZ^2$ exists and $\int_0^t \cZ^1\circ dB = \int_0^t \cZ^2\circ dB$, $\hP$-a.s., for all $t\in [0,T]$, then $\cZ^1=\cZ^2$ as subsets of $\hL^2_\hF([0,T]\times \O, \hR^{d\times m})$.
	\end{lem}

{\it Proof}:
(i) The additivity result (\ref{add1}) is given in \cite[Theorem 2.2]{MgRT}. (ii) Since $\cZ^1,\cZ^2$ are bounded subsets of $\hL^2_\hF([0,T]\times \O,\hR^{d\times m})$, it can be checked that $\cZ^1\ominus\cZ^2$ is also bounded. Moreover, $\cZ^1\ominus\cZ^2$ is convex and closed as a Hukuhara difference. Since $\hL^2_\hF([0,T]\times \O,\hR^{d\times m})$ is reflexive, we may conclude that $\cZ^1\ominus\cZ^2$ is weakly compact. Hence, $\cZ^1\ominus \cZ^2\in \sK_w(\hL^2_{\hF}([0,T]\times\O,\hR^{d\times m}))$. The proof of the 
identity (\ref{sub1}) follows from the additivity of integral as in the proof of Lemma \ref{stochint-add}.

It remains to prove (iii). We first note that by the property of the Hukuhara difference and the assertion (ii), it suffices to show that
$\int_0^t \cZ \circ dB=0$, $\hP\as$, for all $t\in[0,T]$, implies $\cZ=\{0\}$. To see this, we observe that, for a fixed $t\in[0,T]$, the general stochastic integral $\int_0^t \cZ \circ dB=0$, $\hP$-a.s., amounts to saying, by definition, that $S_{\cF_t}^2(\int_0^t\cZ\circ dB)=\overline{\rm dec}_{\cF_t}(\cJ_t[\cZ])=\{0\}$, which is obviously equivalent to $\cJ_t[\cZ]=\{0\}$. In other words, we have $\int_0^t z_sdB_s=0$, $\hP$-a.s., for all $z\in\cZ$. But since this holds for any $t\in[0,T]$, and since the integral $M^z_t:=\int_0^t z_sdB_s$, $t\in[0,T]$, is a continuous martingale, we can conclude that $\hP\{M^z_t=0 \mbox{
for all $t\in[0,T]$}\}=1$ for each $z\in\cZ$. This leads to that $z\equiv 0$, $\hP$-a.s., for all $z\in\cZ$, that is, $\cZ=\{0\}$.
\qed

In the proof of Lemma \ref{genstochint-add}(iii), the existence of the Hukuhara difference $\cZ^1\ominus\cZ^2$ is needed in order to obtain the conclusion $\cZ^1=\cZ^2$, and hence $\cZ^1\ominus\cZ^2=\cb{0}$, by using Lemma \ref{genstochint-add}(ii). To remove this assumption, we will pass to a quotient space of $\sK_w(\hL^2_{\hF}([0,T]\times\O,\hR^{d\times m}))$ in which two sets of processes are considered identical if they yield the same It\^{o} integral. To make this idea precise, let us define a relation $\cong$ on $\sK_w(\hL^2_{\hF}([0,T]\times\O,\hR^{d\times m}))$ by 
\bea
\label{cong}
\cZ^1\cong\cZ^2\quad\Leftrightarrow\quad \int_0^t\cZ^1\circ dB=\int_0^t \cZ^2\circ dB\quad \hP\as\text{ for all }t\in[0,T].
\eea
It is easy to see that $\cong$ is an equivalence relation on $\sK_w(\hL^2_{\hF}([0,T]\times\O,\hR^{d\times m}))$; let us denote $\hK_w(\hL^2_{\hF}([0,T]\times\O,\hR^{d\times m}))$ to be the set of all equivalence classes of $\cong$. For a class $\cZ\in \hK_w(\hL^2_{\hF}([0,T]\times\O,\hR^{d\times m}))$, we define its stochastic integral $\{\int_0^t\cZ\circ dB\}_{t\in[0,T]}$, as the stochastic integral of any member of $\cZ$, which is uniquely defined up to modifications. Hence, for $\cZ^1,\cZ^2\in\hK_w(\hL^2_{\hF}([0,T]\times\O,\hR^{d\times m}))$, if $\int_0^t \cZ^1\circ dB = \int_0^t \cZ^2\circ dB$ $\hP$-a.s., for all $t\in [0,T]$, then $\cZ^1=\cZ^2$ in $\hK_w(\hL^2_\hF([0,T]\times \O, \hR^{d\times m}))$.

For future use, let us extend the definition of Minkowski addition for the new space. For $\cZ,\hat{\cZ}\in \hK_w(\hL^2_\hF([0,T]\times \O, \hR^{d\times m}))$, we define
\bea
\cZ+\hat{\cZ} :=\{\cZ^1+\hat{\cZ}^1\mid \cZ^1\in\cZ,\hat{\cZ}^1\in\hat{\cZ} \},
\eea 
which is well-defined since $\cZ^1+\hat{\cZ}^1\cong \cZ^2+\hat{\cZ}^2$ whenever $\cZ^1,\cZ^2\in\cZ$ and $\hat{\cZ}^1,\hat{\cZ}^2\in\hat{\cZ}$. Then, $\ominus$ has an obvious definition on $\hK_w(\hL^2_\hF([0,T]\times \O, \hR^{d\times m}))$ by \eqref{Hukudiff}. With these definitions, Lemma \ref{genstochint-add} can be rewritten for $\hK_w(\hL^2_\hF([0,T]\times \O, \hR^{d\times m}))$ except that in (iii), the existence of the Hukuhara difference is not needed.

The next corollary is an important observation.

\begin{cor}
	\label{ztimeadd}
	Suppose that $\hP$ is a nonatomic probability measure. Let $\cZ\in\sK_w(\hL^2_{\hF}([0,T]\times\O,\hR^{d\times m}))$ be a nonempty set of processes and $t\in[0,T]$. Then, it holds almost surely that
	\bea\label{ztime1}
	\int_0^T \cZ \circ dB \subset \int_0^t \cZ\circ dB + \int_t^T \cZ\circ  dB,
	\eea
Moreover, if $\cZ$ is decomposable, then it holds almost surely that
	\bea\label{ztime2}
	\int_0^T \cZ \circ dB = \int_0^t \cZ\circ dB + \int_t^T \cZ\circ  dB,\ \int_t^T \cZ \circ dB = \int_0^T \cZ\circ dB \ominus \int_0^t \cZ\circ dB.
	\eea
	\end{cor}

{\it Proof}:
Let $t\in[0,T]$. By \cite[Lemma 3.3.4]{Kis2020}, we have $\cZ\subset\1_{[0,t]}\cZ +\1_{(t,T]}\cZ$, and equality holds when $\cZ$ is decomposable. Applying Lemma \ref{genstochint-add}(i) together with the monotonicity of the integral with respect to $\subset$, the relations in \eqref{ztime1} and \eqref{ztime2} hold. 
\qed

\begin{rem}
\label{remark3.3}
{\rm The essence of  Corollary \ref{ztimeadd} is that, unlike Corollary \ref{dtint-time}, the temporal-additivity $\int_0^T=\int_0^t+\int_t^T$ 
is {\it not} necessarily true in the case of generalized stochastic integrals for lack of decomposability of the integrand. In particular, the Hukuhara difference $\int_0^T \cZ\circ dB\ominus\int_0^t \cZ\circ dB$ may not exist in general. This peculiar feature of generalized stochastic integrals will be particularly felt when we study the set-valued BSDEs in \S\ref{main}.
\qed}
\end{rem}

\section{Some Important Estimates}
\setcounter{equation}{0}

In this section we establish some important estimation regarding set-valued stochastic integrals and  their conditional expectations. These
estimates, albeit conceivable, need justifications given the special natures of the set-valued stochastic analysis, as well as the lack of a vector 
space structure in general. Some of the arguments are following those in \cite{K} closely, but we nevertheless provide the details for the sake of completeness.

Recall the set  $\sK(\hR^d)$, the collection of all nonempty convex compact subsets of $\hR^d$. For $p\in [1,+\infty)$ and $X_1, X_2\in \scL^p_\cF(\O,\sK(\hR^d))$, define
\bea
\label{Hpdist}
\cH_p(X_1,X_2):=(\hE[ h^p(X_1, X_2)])^{\frac{1}{p}}.
\eea 
The following result is a strengthened version of \cite[Theorem 2.4.1]{K} in the $\hL^2$ sense.
%

\begin{lem}
	\label{cond-h2}
	Let $X_1,X_2\in \scL^2_{\cF}(\O,\sK(\hR^d))$, and $\cG \subset \cF$ be a sub-$\sigma$-algebra. Then, one has
	\bea
	\label{ineq0}
	h^2(\hE[X_1|\cG],\hE[X_2|\cG]) \leq \hE[h^2(X_1,X_2)|\cG], \qq \hP\mbox{-a.s.}
	\eea
	In particular, the following inequalities hold:
	\bea
	\label{ineq1}
	\cH_2(\hE[X_1|\cG],\hE[X_2|\cG]) &\leq& \cH_2(X_1,X_2); \ms \\
	\label{ineq2}
	\norm{\hE[X_1|\cG]}^2&\leq& \hE[\norm{X_1}^2|\cG], \qq \hP\mbox{-a.s.}
	\eea 
\end{lem}

{\it Proof.} Let us introduce the notation $\hE[\xi:D]:=\hE[\xi \1_D]$ for $\xi\in\hL^1_\cF(\O,\hR)$ and $D\in\cF$. Note that (\ref{ineq0}) is equivalent to
\bea
\label{contequiv}
\hE\big[ h^2(\hE[X_1|\cG],\hE[X_2|\cG]):D\big] \leq \hE\big[h^2(X_1,X_2):D\big], \qq  D\in\cG.
\eea
Let $D\in\cG$, and define $C:=\{\o\in\O:\bar{h}(\hE[X_1|\cG](\o),\hE[X_2|\cG](\o)) \geq\bar{h}(\hE[X_1|\cG](\o),\hE[X_2|\cG](\o))\}$.
Clearly, by the definition of conditional expectation, $C\in \cG$. Now we can write
\bea
\label{hausdorffdec}
\hE[ h^2(\hE[X_1|\cG],\hE[X_2|\cG]):\neg D]\neg &=& \hE[\bar{h}^2(\hE[X_1|\cG],\hE[X_2|\cG]): D\cap C]\notag \\
&&+\hE[\bar{h}^2(\hE[X_1|\cG],\hE[X_2|\cG]):{D\cap C^c}].
\eea
Repeatedly applying \cite[Theorem 2.3.1]{K}  (see also \cite[Theorem 2.2]{CondE}), we obtain
\bea
\label{Ebarh}
&&\hE\sqb{ \bar{h}^2\of{\hE\sqb{X_1|\cG},\hE\sqb{X_2|\cG}}:{D\cap C}}=\int_{D\cap C} \sup_{x\in\hE\sqb{X_1|\cG}(\o)}d^2(x,\hE\sqb{X_2|\cG}(\o))\hP(d\o)\nonumber\\
&=&\neg\sup_{\eta\in S(\hE\sqb{X_1|\cG})}\hE \sqb{d^2(\eta,\hE\sqb{X_2|\cG})\neg:\neg{D\cap C}}
=\sup_{\eta\in \{\hE\sqb{\varphi|\cG}: \varphi\in S(X_1)\}}\hE \sqb{d^2(\eta,\hE\sqb{X_2|\cG})\neg:\neg{D\cap C}}\nonumber\\
&=&\neg\sup_{\varphi\in S(X_1)}\hE \sqb{d^2(\hE\sqb{\varphi|\cG},\hE\sqb{X_2|\cG})\neg:\neg{D\cap C}}\nonumber\\ 
&=&\sup_{\varphi\in S(X_1)}\int_{D\cap C} \inf_{y\in \hE\sqb{X_2|\cG}(\o) }|\hE[\vf|\cG](\o) -y|^2 \hP(d\o) \nonumber\\
&=&\neg\sup_{\varphi\in S(X_1)}\inf_{\psi\in S(X_2)}\hE\sqb{ \abs{\hE\sqb{\varphi|\cG}-\hE\sqb{\psi|\cG}}^2\neg:\neg{D\cap C}}\nonumber\\
&=&\sup_{\varphi\in S(X_1)}\inf_{\psi\in S(X_2)}\hE\sqb{ \abs{\hE\sqb{\varphi-\psi|\cG}}^2 :{D\cap C}}\nonumber\\
&\leq&\neg \sup_{\varphi\in S(X_1)}\inf_{\psi\in S(X_2)}\hE\sqb{ \hE\sqb{\abs{\varphi-\psi}^2\mid \cG}\neg:\neg{D\cap C}}
= \sup_{\varphi\in S(X_1)}\inf_{\psi\in S(X_2)}\hE\sqb{ \abs{\varphi-\psi}^2:{D\cap C}}\nonumber\\
&=&\hE\sqb{\bar{h}^2(X_1,X_2):{D\cap C}}\leq\hE\sqb{h^2(X_1,X_2):{D\cap C}}.
\eea
Here in the above, the inequality is due to the conditional version of Jensen's inequality. Similarly we also have
$\hE\sqb{ \bar{h}^2\of{\hE\sqb{X_1|\cG},\hE\sqb{X_2|\cG}}:{D\cap C^c}}\leq \hE\sqb{h^2(X_1,X_2):{D\cap C^c}}$. Combining the two inequalities with \eqref{hausdorffdec}, we obtain \eqref{contequiv} and hence \eqref{ineq0}. Then, \eqref{ineq1} is immediate from \eqref{ineq0}. Finally, \eqref{ineq2} follows from \eqref{ineq0} by taking $X_2\equiv\{0\}$.
\qed

Next, we present a H\"older-type of inequality regarding the Aumann integral. A similar inequality appears in \cite[Theorem 2.1]{KisMic16} for a special class of integrands. For completeness, we provide a full proof here for our version.


\begin{prop}
	\label{timeint}
	Let $\Phi^1, \Phi^2\in \scL^{2}_{\hF}([0,T]\times\O,\sK(\hR^d))$, and 
	$t\in [0,T]$. Then, it holds that
	\bea
	\label{Holder}
	h^2\of{\int_t^T\Phi^1_sds,\int_t^T\Phi^2_sds}\leq (T-t)\int_t^T h^2(\Phi^1_s,\Phi^2_s)ds, \qq \hP\as
	\eea
\end{prop}

{\it Proof.} Recalling the definition of the Hausdorff metric $h$,  it suffices to show that
\bea
\label{equiv}
\left\{\ba{lll}
\dis \bar{h}^2\of{\int_t^T\Phi^1_sds,\int_t^T\Phi^2_sds}\leq (T-t)\int_t^T \bar{h}^2(\Phi^1_s,\Phi^2_s)ds,\ms\\
\dis \bar{h}^2\of{\int_t^T\Phi^2_sds,\int_t^T\Phi^1_sds}\leq (T-t)\int_t^T \bar{h}^2(\Phi^2_s,\Phi^1_s)ds, 
\ea\right. \q \hP\as
\eea
By symmetry, we shall check only the first inequality in (\ref{equiv}). To begin with, we first note
that the statement is equivalent to showing, for every $D\in\cF_T$, that 
\bea
\label{est-eq}
\hE\Big[\bar{h}^2\of{\int_t^T\Phi^1_sds,\int_t^T\Phi^2_sds}:D\Big]\leq (T-t)\hE\Big[\int_t^T \bar{h}^2(\Phi^1_s,\Phi^2_s)ds:D\Big].
\eea

To see (\ref{est-eq}), we first note that, similar to (\ref{Ebarh}), we have
\bea
\label{est-eq1}
\hE\Big[ \bar{h}^2\Big(\int_t^T\Phi^1_sds,\int_t^T\Phi^2_sds\Big)\colon D\Big]
=\sup_{\eta_1\in \ol{\dec} J_{t,T}(S^2_{\mathbb{F}}(\Phi^1))}\inf_{\eta_2\in \ol{\dec} J_{t,T}(S^2_{\mathbb{F}}(\Phi^2))}\hE[|\eta_1-\eta_2|^2\colon  D].
\eea
%
%
Next, by the standard H\"older's inequality we have 
\bea
\label{est-eq2}
&&\sup_{\eta_1\in J_{t,T}(S^2_{\mathbb{F}}(\Phi^1))}\inf_{\eta_2\in  J_{t,T}(S^2_{\mathbb{F}}(\Phi^2))}\hE[|\eta_1-\eta_2|^2:D]
\nonumber \\
&&=\sup_{\varphi^1\in S^2_{\mathbb{F}}(\Phi^1)}\inf_{\varphi^2\in  S^2_{\mathbb{F}}(\Phi^2)}\hE\Big[\Big|J_{t,T}(\varphi^1)-J_{t,T}(\varphi^2)\Big|^2:D\Big]\\
&&\le (T-t)\sup_{\varphi^1\in S^2_{\mathbb{F}}(\Phi^1)}\inf_{\varphi^2\in  S^2_{\mathbb{F}}(\Phi^2)}\hE\Big[\int_t^T|\varphi^1_s-\varphi^2_s|^2 ds:D\Big]. \nonumber
\eea		
Now, for given $D\in\cF_T$, we consider the probability space $(D, \cF^D_T, \hP^D)$, where $\cF_T^D:=\cb{C\cap D: C\in \cF_T}$, 
and $\hP^D(C)=[\hP(C)/\hP(D)]\1_{\{\hP(D)>0\}}$, $C\in\cF^D_T$.  We also define the filtration $\hF^D=\{\cF^D_t\}_{t\in[0,T]}$
in a similar way. Applying  \cite[Theorem 2.3.1]{K} again, we have 
%
\bea
\label{est-eq3}
&&\sup_{\varphi^1\in S^2_{\hF}(\Phi^1)}\inf_{\varphi^2\in  S^2_{\hF}(\Phi^2)}\hE\Big[\int_t^T|\varphi^1_s-\varphi^2_s|^2 ds: D \Big]=\sup_{\varphi^1\in S^2_{\hF}(\Phi^1)}\inf_{\varphi^2\in  S^2_{\hF}(\Phi^2)}\hE^{\hP^D}\Big[\int_t^T|\varphi^1_s-\varphi^2_s|^2 ds\Big]\nonumber\\
&=&\sup_{\varphi^1\in S^2_{\hF^D}(\Phi^1)}\inf_{\varphi^2\in  S^2_{\hF^D}(\Phi^2)}\int_{D\times[t,T]}\abs{\varphi^1_s(\o)-\varphi^2_s(\o)}^2\hP^D(d\o) ds\\
&=&\int_{D\times[t,T]}\sup_{x\in \Phi^1_s(\o)}\inf_{y\in  \Phi^2_s(\o)}|x-y|^2\hP^D(d\o)ds=
\int_{D\times[t,T]}\bar{h}^2(\Phi^1_s(\o),\Phi^2_s(\o))\hP^D(d\o) ds\nonumber\\
&=&\hE^{\hP^D}\Big[\int_t^T\bar{h}^2(\Phi^1_s,\Phi^2_s)ds\Big]=\hE\Big[\int_t^T \bar{h}^2(\Phi^1_s,\Phi^2_s)ds:D\Big]. \nonumber
\eea

Let $\a_D:=(T-t)\hE\Big[\int_t^T \bar{h}^2(\Phi^1_s,\Phi^2_s)ds:D\Big]$. Combining (\ref{est-eq2}) and (\ref{est-eq3}), we have
\bea
\label{est-eq4}
\sup_{\eta_1\in J_{t,T}(S^2_{\mathbb{F}}(\Phi^1))}\inf_{\eta_2\in  J_{t,T}(S^2_{\mathbb{F}}(\Phi^2))}\hE[|\eta_1-\eta_2|^2:D]\le \a_D.
\eea

Next, we show that (\ref{est-eq4}) implies that 
\bea
\label{est-eq5}
\sup_{\eta_1\in \dec J_{t,T}(S^2_{\hF}(\Phi^1))}\inf_{\eta_2\in  J_{t,T}(S^2_{\hF}(\Phi^2))}\hE\Big[|\eta_1-\eta_2|^2:D\Big]\leq  \a_D, 
\eea
For any $\eta_1\in \dec J_{t,T}(S^2_{\hF}(\Phi^1))$ we write $\eta_1=\sum_{i=1}^m \1_{D_i}\eta_{1,i}$ for some $D_1,\ldots,D_m\in\cF_T$ partitioning $\O$, and $\eta_{1,1},\ldots,\eta_{1,m}\in J_{t,T}(S^2_{\hF}(\Phi^1))$. Then, for $\eta_2\in J_{t,T}(S^2_{\hF}(\Phi^2))$ we can apply
Jensen's inequality to get
\bea
\label{est-eq6}
\hE[|\eta_1-\eta_2|^2:D]&=&\hE^{\hP^D}\Big[\Big|\sum_{i=1}^m \1_{D_i}(\eta_{1,i}-\eta_2)\Big|^2\Big]\le  \hE^{\hP^D}\Big[\sum_{i=1}^m \1_{D_i}|\eta_{1,i}-\eta_2|^2\Big]\nonumber\\\
&=&\sum_{i=1}^m \hE\big[\1_{D\cap D_i}|\eta_{1,i}-\eta_2|^2\big].
\eea
Since $\eta_1$ and $\eta_2$ are arbitrary, we deduce from (\ref{est-eq6}) and (\ref{est-eq4}) that 
\beaa
&&\sup_{\eta_1\in \dec J_{t,T}(S^2_{\hF}(\Phi^1))}\inf_{\eta_2\in J_{t,T}(S^2_{\hF}(\Phi^2))}\hE[|\eta_1-\eta_2|^2:D]\\
&\le &\sum_{i=1}^m\sup_{\eta_{1,i}\in  J_{t,T}(S^2_{\hF}(\Phi^1))}\inf_{\eta_2\in J_{t,T}(S^2_{\hF}(\Phi^2))} \hE\big[|\eta_{1,i}-\eta_2|^2:{D\cap D_i}\big]\le \sum_{i=1}^m \a_{D\cap D_i}=\a.
\eeaa
This proves \eqref{est-eq5}.
Noting that  $\ol{\dec} J_{t,T}(S^2_{\hF}(\Phi^2))\supset J_{t,T}(S^2_{\hF}(\Phi^2))$,   (\ref{est-eq5}) implies that
\bea
\label{est-eq51}
\sup_{\eta_1\in \dec J_{t,T}(S^2_{\hF}(\Phi^1))}\inf_{\eta_2\in  \ol{\dec} J_{t,T}(S^2_{\hF}(\Phi^2))}\hE\Big[|\eta_1-\eta_2|^2:D\Big]\leq  \a_D.
\eea
Finally, we claim that \eqref{est-eq51} implies 
\bea
\label{est-eq52}
\sup_{\eta_1\in \ol{\dec} J_{t,T}(S^2_{\hF}(\Phi^1))}\inf_{\eta_2\in  \ol{\dec} J_{t,T}(S^2_{\hF}(\Phi^2))}\hE\Big[|\eta_1-\eta_2|^2:D\Big]\leq  \a_D,
\eea
which, together with (\ref{est-eq1}), would lead to (\ref{est-eq}). Indeed, let $\eta_1\in \ol{\dec} J_{t,T}(S^2_{\hF}(\Phi^1))$, and 
let  $\{\eta_1^n\}_{n\in\hN}\subset \dec J_{t,T}(S^2_{\hF}(\Phi^1))$ be a sequence that converges to $\eta_1$ (strongly) in $\hL^2_{\cF_T}(\O,\hR^d)$. Let $\varepsilon>0$.  For each $n\in\hN$, thanks to \eqref{est-eq51}, we may find $\eta_2^n\in \ol{\dec} J_{t,T}(S^2_{\hF}(\Phi^2))$ such that
\bea 
\label{est-eq7}
\hE[|\eta^n_1-\eta^n_2|^2:D]<\a_D+\varepsilon.
\eea 
By Remark \ref{propint}, $\{\eta_2^n\}_{n\in\hN}$ is a bounded sequence in $\hL^2_{\cF_T}(\O,\hR^d)$; hence, by Banach-Saks theorem, it has a subsequence $\{\eta_2^{n_k}\}_{k\in\hN}$ for which the sequence $\{\ol{\eta}_2^k\}_{k\in\hN}$ converges to some $\ol{\eta}_2\in \hL^2_{\cF_T}(\O,\hR^d)$ strongly, where $\ol{\eta}_2^k:=\frac1k \sum_{\ell=1}^k\eta_2^{n_{\ell}}$ is the Ces\`{a}ro average, for $k\in\hN$. Moreover, since $\ol{\dec} J_{t,T}(S^2_{\hF}(\Phi^2))$ is a closed convex set, all Ces\`{a}ro averages and their limit $\ol{\eta}_2$ belong to $\ol{\dec} J_{t,T}(S^2_{\hF}(\Phi^2))$. The strong convergence of $\{\eta_1^n\}_{n\in\hN}$ implies that $\{\ol{\eta}_1^k\}_{k\in\hN}\subset \ol{\dec} J_{t,T}(S^2_{\hF}(\Phi^1))$ converges to $\eta_1$ strongly in $\hL^2_{\cF_T}(\O,\hR^d)$, where $\ol{\eta}_1^k:=\frac1k \sum_{\ell=1}^k\eta_1^{n_{\ell}}$, $k\in\hN$. By \eqref{est-eq7}, we have
\beaa
\label{est-eq8}
\hE[|\ol{\eta}_1^k-\ol{\eta}_2^k|^2:D]
\le \Big(\frac1k\sum_{\ell=1}^k \Big(\hE[|\eta_1^{n_\ell}-\eta_2^{n_\ell}|^2:D]\Big)^{\frac12}\Big)^2<\a_D+\varepsilon,\q k\in\hN.
\eeaa 
Thus,
\beaa 
\label{est-eq9}
(\hE[|\eta_1-\ol{\eta}_2|^2:D])^{\frac12}&\le & (\hE[|\eta_1-\ol{\eta}^k_1|^2:D])^{\frac12}+(\hE[|\ol{\eta}^k_1-\ol{\eta}^k_2|^2:D])^{\frac12}+(\hE[|\ol{\eta}^k_2-\ol{\eta}_2|^2:D])^{\frac12}\\
&\le & (\hE[|\eta_1-\ol{\eta}^k_1|^2])^{\frac12}+(\a_D+\varepsilon)^{\frac12}+(\hE[|\ol{\eta}^k_2-\ol{\eta}_2|^2])^{\frac12},
\eeaa 
and letting $k\rightarrow\infty$ yields
\beaa
\inf_{\eta_2\in  \ol{\dec} J_{t,T}(S^2_{\hF}(\Phi^2))}\hE[|\eta_1-\eta_2|^2:D]\le \hE[|\eta_1-\ol{\eta}_2|^2:D]\leq \a+\varepsilon.
\eeaa
Since $\varepsilon>0$ and $\eta_1\in \ol{\dec} J_{t,T}(S_{\hF}(\Phi^1))$ are arbitrary, \eqref{est-eq52} follows, concluding the proof.
\qed


\section{Set-Valued Martingales and their Integral Repsentations}

Using the notion of {\it conditional expectation} in \S\ref{svrv}, one can define {\it set-valued martingales} as follows.
%
We say that a set-valued process $M=\{M_t\}_{t\in[0,T]}$ is a {\it set-valued $\hF$-martingale} if $M \in \scL^0_\hF([0,T] \times \Omega, \sC(\hR^d))$, $M_t\in \sA^1_{\cF_t}(\O,\sC(\hR^d))$, and
$M_{s} = \hE[M_t|\cF_s]$ for all $0\le s\le t$.  $M$ is called {\it square-integrable} if $M_t \in\sA^2_{\cF_t}(\O,\sC(\hR^d)) $, $0\le t \le T$, and {\it uniformly square-integrably bounded} if there exists $\ell\in\hL^2(\O,\hR_+)$ such that $\sup_{t\in[0,T]}\norm{M_t(\cd)}\leq \ell(\cd)$ a.s.

We note that, if $M$ is a  square-integrable  set-valued martingale, then for each $t\in[0,T]$, the set  of  square-integrable selectors, $S^2_{\cF_t}(M_t)$, is decomposable. On the other hand, we consider the set of all square-integrable martingale selectors, that is, all $d$-dimensional $\hF$-martingales $f=\{f_t\}_{t\in[0,T]}$ such that $f_t\in S^2_{\cF_t}(M_t)$, $t\in[0,T]$, and denote it by $MS(M)$.  If $M$ is convex-valued, then it is known that $MS(M)\neq\emptyset$ (cf. \cite[\S3]{MgRT}).  For $t\in[0,T]$, consider the $t$-section of $MS(M)$, defined as
$P_t[MS(M)]:=\{f_t: f\in MS(M)\}\subset \hL^2_{\cF_t}(\O, \hR^d)$. We remark that
the two sets $S^2_{\cF_t}(M_t)$
(the selectors of the $t$-section) and $P_t[MS(M)]$ (the $t$-section of the selectors) are quite different. In particular, the former is
known to be decomposable, but the latter is not. However, the following relation holds (see \cite[Proposition~3.1]{MgRT}):
\bea
\label{Projdecom}
S^2_{\cF_t}(M_t)=\overline{\dec}_{\cF_t}( P_t[MS(M)]), \qq t\in [0,T],
\eea
where $\overline{\dec}_{\cF_t}$ denotes the closed decomposable hull with respect to $\hL^2_{\cF_t}(\O,\hR^d)$. 

%

\subsection{Representation of Martingales with Trivial Initial Value}\label{sec:zero_init}

In what follows we assume that $\hF=\hF^B$, for some $\hR^m$-valued Brownian motion $B=\{B_t\}_{t\in[0,T]}$. The fundamental building block of the theory  of Backward SDE is the celebrated {\it Martingale Representation Theorem}, which states that
 every square-integrable $\hF$-martingale can be written, uniquely, as a stochastic integral against $B$, whence continuous. 
There is a similar result for set-valued martingales (see \S\ref{svsp}), which we now describe. 

%


Let $M$ be a convex-valued set-valued $\hF$-martingale that is square-integrable, i.e., $M_t\in \sA^2_{\cF_t}(\O,\sC(\hR^d))$ for each $t\in [0,T]$.
Then for each $y\in MS(M)$, by standard martingale representation theorem, there exists unique $z^y\in \hL^{2}_\hF([0,T],\hR^{d\times m})$, such that 
$y_t=\int_0^t z^y_sdB_s$, $t\in[0,T]$, $\hP$-a.s. Denote $\cZ^M:=\{z^y: y\in MS(M)\}\in \sP(\hL^2_\hF([0,T],\hR^{d\times m}))$. 
%
\begin{rem}\label{mtgrem}
	{\rm
	We should note that while  a set-valued martingale always gives rise to a set of vector-valued martingales, i.e., stochastic integrals, 
	not every 
	set of vector-valued martingales can be realized as $MS(M)$ for some set-valued martingale $M$.
	\qed}
	\end{rem}

The following {\it Set-valued Martingale Representation Theorem} is due to \cite{MgRT}.

\begin{thm}[Kisielewicz {\cite[Proposition 4.1, Theorem 4.2]{MgRT}}]
\label{SVMgRT}
For every convex-valued square-integrable set-valued martingale $M = \{M_t\}_{t\in [0, T]}$ with $M_0=\{0\}$, there exists $\cZ^M\in \sP(\hL^2_\hF([0,T]\times\O,\hR^{d\times m}))$ such that $M_t = \int_0^t \cZ^M  \circ dB$, $\hP$-a.s. $t\in [0, T]$.  If $M$ is also uniformly square-integrably bounded, then $\cZ^M$ is a convex weakly compact set, that is, $\cZ^M\in \sK_w(\hL^2_{\hF}([0,T]\times\O,\hR^{d\times m}))$.
\end{thm}
\begin{rem}
\label{remMgRT}
{\rm
(i) We first note that in the set-valued martingale representation,  the ``martingale integrand" $\cZ^M$ may not be a
measurable set-valued process. In fact, if the set-valued martingale is square-integrably bounded, then the integrand $\cZ^M$ {\it cannot} be decomposable unless it is a singleton (see \cite[Corollary 5.3.2]{Kis2020}). Thus the stochastic integral can only be in the generalized sense.  But on the other hand, if $\cZ^M$ is not decomposable, 
then the temporal-additivily of the set-valued stochastic integral  fails in general (see, Corollary \ref{ztimeadd}). Such a conflict leads to some fundamental difficulties for the study of  set-valued BSDEs, and it does not seem to be amendable unless some more general framework of set-valued stochastic integrals is established. 


(ii) If  $\O$ is separable, then there exists a sequence $\{z_n\}_{n\geq 1} \subset \hL^2_\hF([0.T],\hR^{d
\times m})$ such that $M_t=\cl\{\int_0^t z^n_s dB_s\}_{n\ge 1}$, $S^2_{\cF_t}(M_t)=\overline{\dec}_{\cF_t}\{\int_0^t z^n_s dB_s\}_{n\ge 1}$, $t\in[0,T]$ (see \cite[Theorem 4.3]{MgRT}). 

(iii) If $M$ is a uniformly square-integrably bounded martingale and $\hP$ is nonatomic, then there exists a sequence $\{z_n\}_{n\geq 1} \subset \hL^2_\hF([0.T],\hR^{d
	\times m})$ such that $M_t=\overline{\co}\{\int_0^t z^n_s dB_s\}_{n\ge 1}$ for all $t\in[0,T]$ (see \cite[Theorem 4.3]{MgRT}). 

(iv) In light of the equivalence relation $\cong$ in \eqref{cong}, in the last part of Theorem \ref{SVMgRT}, we can easily conclude that such $\cZ^M$ is unique in $\hK_w(\hL^2_{\hF}([0,T]\times\O,\hR^{d\times m}))$. Indeed, if there exist $\cZ^M_1$ and $\cZ^M_2$ in $\sK_w(\hL^2_{\hF}([0,T]\times\O,\hR^{d\times m}))$ such that $\int_0^t \cZ^M_1\circ dB =\int_0^t \cZ^M_2\circ dB=M_t$, 
$t\in[0,T]$, then $\cZ^M_1\cong\cZ^M_2$, that is, they correspond to the same element of $\hK_w(\hL^2_{\hF}([0,T]\times\O,\hR^{d\times m}))$ and we may denote this element by $\cZ^M$ with a slight abuse of notation.

(v) Unlike usual stochastic integrals, set-valued stochastic integrals do not always generate set-valued martingales. In fact, given a nonempty set $\cZ\in \sK_w(\hL^2_{\hF}([0,T]\times\O,\hR^{d\times m}))$ (or $\cZ\in \hK_w(\hL^2_{\hF}([0,T]\times\O,\hR^{d\times m}))$)
of processes, the set-valued process $\{\int_0^t\cZ\circ dB\}_{t\in [0,T]}$ forms a set-valued {\it submartingale} in the sense that $\int_0^u\cZ\circ dB\subset \hE[\int_0^t\cZ\circ dB_s\vert \cF_u]$ for every $0\leq u\leq t\leq T$ (see \cite[Theorem 4.2]{KisMic17}). Nevertheless, the stochastic integrals that appear in Theorem \ref{SVMgRT} are naturally martingales. 
\qed }
\end{rem}

\subsection{Representation of Martingales with General Initial Value}\label{definiteint}

We would like to point out that in Theorem \ref{SVMgRT} it is assumed that $M_0=\{0\}$. Such a seemingly benign assumption actually has some severe consequences. In particular, as it was pointed out recently in \cite{ZhangYano}, a set-valued martingale whose initial value is a singleton is essentially a vector-valued martingale. Therefore, Theorem \ref{SVMgRT} actually is not a suitable tool for 
the study of set-valued BSDEs with non-singleton terminal values. The main purpose of this subsection is to establish a refined version of set-valued martingale representation theorem for 
set-valued martingales with general (non-singleton) initial values. 

Our idea is to extend the notion of Aumman-It\^o integral so that 
it is a martingale but its expectation is not necessarily zero (see \cite[Example 3.1]{ZhangYano} for the set-valued delemma). To this end, for any $t\in [0,T]$, we consider the space
$\hR_t:=\hL^2_{\cF_t}(\O,\hR^d)\times \hL_{\hF}^2([t,T]\times\O,\hR^{d\times m})$.

 Given a process $z=\{z_u\}_{u\in [0,T]}$, we denote $z^{t,T}:=(z_u)_{u\in[t,T]}$ to be the restriction of $z$ onto the interval $[t,T]$, and
 define a mapping $F^t\colon\hR_0\mapsto\hR_t$ by
\[
F^t(x,z):=\Big(x+\int_0^t z_s dB_s, z^{t,T}\Big),\quad (x,z)\in\hR_0.
\]
We have the following result.
\begin{lem}\label{lemFJ}
For given $t\in[0,T]$ and $(\xi,z^t)\in \hR_t$, define a process  $\cJ^t(\xi,z^t)=\{\cJ^t_u(\xi,z^t)\}_{u\in[0,T]}$:
\bea\label{J-int}
\cJ^t_u(\xi,z^t):=\hE[\xi|\cF_u]\1_{[0,t)}(u)+\of{\xi+\int_t^u z^t_s dB_s}\1_{[t,T]}(u),\quad u\in [0,T].
\eea
Then, $\cJ^t(\xi, z^t)$ is an $\hF$-martingale on $[0,T]$. Moreover, it holds that 
 $\cJ^t\circ F^t = \cJ^0$ on $\hR_0$.
\end{lem}

{\it Proof.} That $\cJ^t(\xi, z^t)$ is a martingale is obvious. To check the identity, let $(x,z)\in \hR_0$. Following the definitions of $\cJ^t, F^t, \cJ^0$, we have
\beaa
\cJ^t_u(F^t(x,z))&=&\cJ^t_u\Big(x+\int_0^t z_sdB_s, z^{t,T}\Big)\\
&=&\hE\Big[x+\int_0^t z_sdB_s\Big| \cF_u\Big]\1_{[0,t)}(u)+\Big(x+\int_0^t z_s dB_s + \int_t^u z_sdB_s\Big)\1_{[t,T]}(u)\\
&=&\Big(x+\int_0^u z_sdB_s\Big)\1_{[0,t)}(u)+\Big(x+\int_0^u z_sdB_s\Big)\1_{[t,T]}(u)=\cJ^0_u(x,z),
\eeaa 
for every $u\in[0,T]$. Hence, $\cJ^t(F^t(x,z))=\cJ^0(x,z)$.
\qed 

Next, let $\cR\subset \hR_0$ be a nonempty set and $t\in[0,T]$. By virtue of Theorem \ref{decomp}, there exists 
a set-valued random variable in $\hL^2_{\cF_t}(\O,\sC(\hR^d))$, denoted by $\int_{0-}^t\cR\circ dB$, such that
\bea
\label{stochint-defn}
S^2_{\cF_t}\Big(\int_{0-}^t \cR\circ dB\Big)=\overline{\dec}_{\cF_t}(\cJ^0_t[\cR]).
\eea
We call $\int_{0-}^t\cR\circ dB$ the {\it stochastic integral} of $\cR$. Clearly, such a stochastic integral is an extended version of the generalized Aumann-It\^{o} stochastic integral $\int_0^t \cZ\circ dB$ defined by \eqref{genint}, 
and in particular, the integrand $\cR$ consists of pairs $(x,z)$, which keeps track of the initial values $x$ of the martingales in $\cJ^0[\cR]$, motivating the choice of the notation $\int_{0-}^t$.

To see how the integral $\int_{0-}^t \cR\circ dB$ (or more precisely, $\cR$) can be defined through a set-valued martingale, let
$M=\{M_u\}_{u\in[0,T]}$ be a convex uniformly square-integrably bounded set-valued martingale with respect to $\hF=\hF^B$, 
and $M_0$ is a non-singleton convex set. Let $MS(M)$ be the set of all $\hL^2$-martingale selectors of $M$. By standard martingale representation theorem, for fixed $t\in[0,T]$, each $y\in MS(M)$ can be written as
$y=\cJ^t(\xi,z)$ for a unique pair $(\xi,z)\in\hR_t$. We shall define, for each $t\in[0,T]$,
\bea\label{RMdefn}
\cR_t^M:=\cb{(\xi,z)\in \hR_t\colon \cJ^t(\xi,z)\in MS(M)}; \q\mbox{\rm and}\q \cR^M:=\cR^M_0.
\eea

In what follows,  for  $(\xi,z)\in\hR_t$, we  write $\pi_\xi(\xi,z):=\xi$ and $\pi_z(\xi,z):=z$. (For  convenience, we suppress the dependence of the mappings $\pi_\xi,\pi_z$ on $t$.) Also, if $y=\cJ^t(\xi, z)\in \cR_t^M$, we denote $\pi_\xi(y)=\xi$, 
and $\pi_x(y)=x$, respectively. Furthermore, 
we define $\cZ^M_t:=\pi_z[\cR^M_t]$. The following theorem collect collects various forms of ``time-consistency" properties of the collection $\{\cR_t^M\}_{t\in[0,T]}$, which will be useful in our future discussion.
\begin{prop}[Time-consistency]
\label{integrand-flow}
	 Let $t_1\in[0,T]$. Then, it holds
	\bea \label{tcmain}
	F^{t_1}[\cR^M_{0}]=\cR^M_{t_1}.
	\eea 
	Furthermore, the following relations hold for every $t_2\in (t_1,T]$:
	
	(i) $\cJ^{t_1}[\cR^M_{t_1}]=\cJ^{t_2}[\cR^M_{t_2}]=MS(M)$.
	
	(ii) $\pi_{\xi}[\cR^M_{0}]=M_0$.
	
	(iii) $\pi_\xi[\cR^{M}_{t_1}]=J^{t_1}_{t_1}[\cR^M_{t_1}]=J^{t_2}_{t_1}[\cR^M_{t_2}]=P_{t_1}[MS(M)]$.
	
	(iv) $\pi_\xi[\cR^M_{t_1}]=\{\hE[\xi\;\vert\;\cF_{t_1}]\colon \xi\in \pi_\xi[\cR^M_{t_2}]\}$.
	
	(v) $\cZ^{M}_{t_1}\1_{[t_2,T]}=\cZ^{M}_{t_2}\1_{[t_2,T]}=\cZ^M_0\1_{[t_2,T]}$.
\end{prop}

{\it Proof.} We first prove \eqref{tcmain}. Fix $t_1\in[0,T]$ and let $(x,z)\in \cR^M_0$. By Lemma \ref{lemFJ} and the definition of $\cR^M_0$, we have 
$\cJ^{t_1}(F^{t_1}(x,z))=\cJ^0(x,z)\in MS(M)$.
On the other hand, since $(x,z)\in\hR_0$, we have $F^{t_1}(x,z)\in \hR_{t_1}$,  which implies $F^{t_1}(x,z)\in \cR^M_{t_1}$. Namely, $F^{t_1}[\cR^M_0]\subset \cR^M_{t_1}$. 

Conversely, let $(\hat{\xi},\hat{z})\in \cR^M_{t_1}$. Define a martingale $y_s:=\hE[\hat{\xi}|\cF_s]$, $s\in[0,t_1]$.
By martingale representation theorem, there exists a unique pair $(x,\bar{z})\in\hR_0$ such that 
\[
y_u=x + \int_{0}^u \bar{z}_sdB_s,\quad u\in[0,t_1].
\]
Let $z:=\bar{z}\1_{[0,t_1)}+\hat{z}\1_{[t_1,T]}\in\hL^2_{\hF}([0,T]\times\O,\hR^{d\times m})$. Then, $(x,z)\in\hR_0$ and for every $u\in[0,T]$,
\beaa
\cJ^{0}_u(x,z)&=&x+\int_0^{u}z_sdB_s=\of{x+\int_0^u \bar{z}_sdB_s}\1_{[0,t_1)}(u)+\of{\hat{\xi}+\int_{t_1}^u\hat{z}_sdB_s}\1_{[t_1,T]}(u)\\
&=&\hE[\hat{\xi}\;\vert\;\cF_u]\1_{[0,t_1)}(u)+\of{\hat{\xi}+\int_{t_1}^u\hat{z}_sdB_s}\1_{[t_1,T]}(u)=\cJ_u^{t_1}(\hat{\xi},\hat{z}).
\eeaa
Hence, $\cJ^0(x,z)=\cJ^{t_1}(\hat{\xi},\hat{z})\in MS(M)$, that is, $(x,z)\in \cR^M_0$. Finally,
\[
F^{t_1}(x,z)=\of{x+\int_0^{t_1}z_sdB_s,z^{t_1,T}}=(\hat{\xi},\hat{z}).
\]
So $(\hat{\xi},\hat{z})\in F^{t_1}[\cR^M_0]$. Consequently, we have $\cR^M_{t_1}\subset F^0[\cR^M_0]$, proving \eqref{tcmain}.

We now turn to properties (i)--(v). The proof of (i) is immediate since $\cJ^{t_i}[\cR^M_{t_i}]=MS(M)$ by the definition of $\cR^M_{t_i}$ for $i\in\{1,2\}$.

To see (ii), let $(x,z)\in\cR^M_0$. Since $\cJ^0(x,z)\in MS(M)$, we have $\pi_\xi(x,z)=x=\cJ^0_0(x,z)\in M_0$. Conversely, since $M$ is a set-valued martingale, $M_0=\hE[M_T|\cF_0]=\hE[M_T]$, thanks to Blumenthal 0-1 law. Hence, by the definition of set-valued expectation, for any $x\in M_0$, there exists $\xi\in S^2_{\cF_T}(M_T)$ such that $x=\hE[\xi]$. Furtherm, by the martingale representation theorem, there exists $z\in\hL^2_{\hF}([0,T]\times\O,\hR^{d\times m})$ such that
\[
\hE[\xi|\cF_u]=x+\int_0^u z_sdB_s = \cJ^0_u(x,z),\quad u\in [0,T]. 
\]
Note that $M$ is a set-valued martingale, we have $\hE[\xi\vert\cF_u]\in S^2_{\cF_u}(M_u)$, $u\in [0,T]$. Hence, $J^0(x,z)=\{\hE[\xi\vert\cF_u]\}_{u\in[0,T]}\in MS(M)$, that is, $(x,z)\in \cR^M_0$, or $x\in\pi_\xi[\cR^M_0]$, proving (ii).

To prove (iii), first note that, for every $(x,z)\in\hR_0$,
\[
\pi_\xi(F^{t_1}(x,z))=\pi_\xi\Big(x+\int_0^{t_1}z_sdB_s,z^{t_1,T}\Big)=x+\int_0^{t_1}z_sdB_s=\cJ^0_{t_1}(x,z).
\]
This implies that
\beaa
\pi_\xi[\cR^M_{t_1}]=\pi_\xi[F^t[\cR^M_0]]=\cb{\cJ^0_{t_1}(x,z)\colon (x,z)\in \cR^M_0}= \cJ^{0}_{t_1}[\cR^M_0]=\cJ^{t_1}_{t_1}\circ F^{t_1}[\cR^M_0]=\cJ^{t_1}_{t_1}[\cR^M_{t_1}],
\eeaa
where the first and last equalities are by \eqref{tcmain} and the fourth equality is due to Lemma \ref{lemFJ}. On the other hand, by the definitions of $P_{t_1}, \cJ^{t_1}$, we see that $P_{t_1}\circ \cJ^{t_1}=\cJ^{t_1}_{t_1}$. Therefore, 
\[
\pi_\xi[\cR^M_{t_1}]=\cJ^{t_1}_{t_1}[\cR^M_{t_1}] = P_{t_1}[\cJ^{t_1}[\cR^M_{t_1}]]=P_{t_1}[\cJ^{t_2}[\cR^M_{t_2}]]=P_{t_1}[MS(M)],
\]
thanks to (i), which concludes the proof of (iii).

To prove (iv), first note that $\hE[P_{t_2}(y)| \cF_{t_1}]=P_{t_1}(y)$ whenever $y=\{y_u\}_{u\in[0,T]}$ is a martingale. Hence, applying (iii) twice, we obtain
\[
\{\hE[\xi|\cF_{t_1}]\colon \xi\in \pi_\xi[\cR^M_{t_2}]\}=\{\hE[\xi|\cF_{t_1}]\colon \xi\in P_{t_2}[MS(M)]\}=P_{t_1}[MS(M)]=\pi_\xi[\cR^M_{t_1}].
\]

Finally, to prove (v), note that, for every $(x,z)\in\hR_0$,
\[
\pi_z(F^{t_1}(x,z))=\pi_z\Big(x+\int_0^{t_1}z_sdB_s,z^{t_1,T}\Big)=z^{t_1,T}.
\]
Hence,
\beaa
\cR^M_{t_1}\1_{[t_2,T]}&=&\pi_z[\cR^M_{t_1}]\1_{[t_2,T]}=\pi_z[F^{t_1}[\cR^M_0]]\1_{[t_2,T]}\\
&=&\cb{z^{t_1,T}\colon (x,z)\in \cR^M_0}\1_{[t_2,T]}=\cb{z\1_{[t_2,T]}\colon z\in \cZ^M_0}=\cZ^M_0\1_{[t_2,T]}.
\eeaa
Taking $t_1=t_2$ above, we also obtain $\cZ^M_{t_2}\1_{[t_2,T]}=\cZ^M_0\1_{[t_2,T]}$.
\qed 

The following theorem is a martingale representation theorem for set-valued martingales with possibly nontrivial initial values, i.e., $M_0$ is a non-singleton deterministic set. 

\begin{thm}\label{MRT-initial}
	Let $M=\{M_u\}_{u\in[0,T]}$ be a convex uniformly square-integrably bounded set-valued martingale with respect to $\hF=\hF^B$. Then, for each $u\in[0,T]$, it holds
	\[
	M_u = \int_{0-}^u \cR^M \circ dB\quad \text{a.s.}
	\]
	Moreover, for each $t\in [0,u]$, it holds that $S^2_{\cF_u}(M_u)=\overline{\dec}_{\cF_u}(\cJ_u^t[\cR^M_t])$.
\end{thm}

{\it Proof}: By Lemma \ref{integrand-flow}(ii),  we have
$\cJ^0_u[\cR^M]=P_u[\cJ^0[\cR^M]]=P_u[MS(M)]$, $u\in[0,T]$.
On the other hand, by \cite[Proposition 3.1]{MgRT}, we have
$\overline{\dec}_{\cF_u}(P_u[MS(M)])=S^2_{\cF_u}(M_u)$.
Combining these with the definition of stochastic integral in \eqref{stochint-defn}, we get
\[
S^2_{\cF_u}\Big(\int_{0-}^u \cR^M\circ dB\Big)=\overline{\dec}_{\cF_u}(\cJ^0_u[\cR^M])=\overline{\dec}_{\cF_u}(P_u[MS(M)])=S^2_{\cF_u}(M_u).
\]
This shows that $M_u=\int_{0-}^u \cR^M\circ dB$ almost surely. The second part of the proposition is an immediate consequence of Lemma \ref{lemFJ}.
\qed

\begin{rem}\label{comparison}
{\rm
	It is interesting to note the relationship between the new stochastic integral $\int_{0-}^u \cR^M\circ dB$  \eqref{stochint-defn} and the generalized Aumann-It\^{o} stochastic integral $\int_0^u \cZ^M \circ dB$  \eqref{genint}, where $\cZ^M:=\cZ^M_0$. Recalling \eqref{JsJ} and \eqref{J-int}, we have
	\beaa
	\cJ^0_u(x,z)=x+\cJ_u(z)\in M_0+\cJ_u[\cZ^M]
	\eeaa 
	for every $(x,z)\in\cR^M$. Hence, $\cJ^0_u[\cR^M]\subset M_0+\cJ_u[\cZ^M]$. After taking closed decomposable hulls, it follows that
	\[
	S^2_u\Big(\int_{0-}^u \cR^M\circ dB\Big)=\overline{\dec}_{\cF_u}(\cJ^0_u[\cR^M])\subset M_0+\overline{\dec}_{\cF_u}(\cJ_u[\cZ^M])=M_0+S^2_{\cF_u}\Big(\int_0^u \cZ^M\circ dB\Big).
	\]
Therefore,
	\[
	\int_{0-}^u \cR^M\circ dB \subset M_0+\int_0^u \cZ^M\circ dB\text{ a.s.}
	\]
	and the reverse inclusion fails to hold in general. When $M_0=\{0\}$, we have $\int_{0-}^u \cR^M \circ dB =\int_0^u \cZ^M\circ dB$ since $\cR^M=\{0\}\times\cZ^M$ in this case.
\qed}
\end{rem}

\begin{rem}\label{comparison2}
	\rm 
	In view of Remark \ref{comparison} and Theorem \ref{MRT-initial},  the new integral $\int_{0-}^u \cR^M\circ dB$ is a non-trivial and necessary extension of the Aumann-It\^{o} stochastic integral, which can be used for the integral representation of any truly set-valued martingale $M$ with a non-zero (non-singleton) initial value $M_0$.
	\qed
\end{rem}

\section{Set-Valued BSDEs}\label{main}
\setcounter{equation}{0}

We are now ready to study the set-valued BSDEs. Assume from now on that $(\O, \cF, \hP, \hF)$ is a filtered probability space on which
is defined an $m$-dimensional standard  Brownian motion $B=\{B_t\}_{t\in[0,T]}$. We assume further that $\hF=\hF^B$, the natural filtration generated by  $B$, augmented by all the $\hP$-null sets  of $\cF$  so that it satisfies the {\it usual hypotheses}. In particular, we may assume
without loss of generality that  $(\O, \cF)=(\hC([0,T]), \sB(\hC([0,T])))$ is the canonical space with $\cF_t=\si\{\o(\cd\wedge t), \o\in\O\}$, $t\in[0,T]$, and $\hP$ is the Wiener measure on $(\O, \cF)$. Hence, $\O$ is separable and $\hP$ is nonatomic.

\subsection{Set-Valued BSDEs in Conditional Expectation Form}\label{sec:condexp}

In this section, we shall focus on the following simplest form of set-valued BSDE:
\bea
\label{BSDEc}
Y_t = \hE\Big[\xi + \int^T_t  f(s, Y_s)ds\ \Big|\ \cF_t\Big], \qq t\in[0,T].
\eea
where $\xi\in \scL^2_{\cF_T}(\O,\sK(\hR^d))$, $f\colon [0,T]\times\O\times \sK(\hR^d)\to \sK(\hR^d)$ is a  set-valued function to be specified later. 
We first give the definition of the solution to the set-valued BSDE (\ref{BSDEc}). 
\begin{defn}
	\label{svbsde0}
	A set-valued process $Y\in \scL^{2}_{\hF}([0,T]\times\O,\sK(\hR^d))$ is called an adapted solution to the set-valued BSDE \eqref{BSDEc}
	if
	\bea
	Y_t = \hE\Big[\xi + \int^T_t  f(s, Y_s)ds\ \Big|\ \cF_t\Big],\qq \hP\as, t\in[0,T].\nonumber 
	\eea
\end{defn}

We shall make use of the the following assumptions on the coefficient $f$.
\begin{assum}\label{assump2}
	The function $f:[0,T]\times\O\times \sK(\hR^d)\to \sK(\hR^d)$ enjoys the following properties:
	
	(i) for fixed $A\in\sK(\hR^d)$, $f(\cd, \cd, A)\in\scL^0_\hF([0,T]\times\O,\sK(\hR^d))$;
	
	(ii)  $f(\cd, \cd, \{0\})\in \scL^{2}_\hF([0,T]\times\O,\sK(\hR^d))$, that is,  
	\bea
	\label{Int}
	\hE\Big[\int_0^T \|f(t, \{0\})\|^2dt\Big]=\hE\Big[\int_0^T h^2(f(t, \{0\}), \{0\})dt\Big]<\infty;
	\eea
	(iii) 
	%
	$f(t,\o,\cdot)$ is Lipschitz, uniformly in $(t,\o)\in[0,T]\times\O$, in the following sense: there exists $K>0$ such that
	\bea
	\label{Lip}
	h(f(t,\o, A),f(t, \o, B))\le Kh(A, B), \qq A, B\in\sK(\hR^d),\ (t,\o)\in [0,T]\times\O.
	\eea
\end{assum}

\begin{rem}\label{Cara}
	{\rm Note that a multifunction $f$ satisfying Assumption \ref{assump2} must be a Carath\'eodory multifunction (see Section \ref{svmappings}), which requires only continuity on the spatial variable. 
		\qed}
\end{rem}

\begin{rem}\label{assumpcon}
{\rm	By Assumption \ref{assump2}, it is easy to check that $\{f(t,Y_t)\}_{t\in[0,T]}\in \scL^{2}_{\hF}([0,T]\times\O,\sK(\hR^d))$ whenever $\{Y_t\}_{t\in[0,T]}\in \scL^{2}_{\hF}([0,T]\times\O,\sK(\hR^d))$.
\qed}
\end{rem}

We shall consider the following standard Picard iteration. Let  $Y^{(0)} \equiv \{0\}$ and for $n\ge 1$, we define $Y^{(n)}$ recursively by
	\bea \label{picardseq}
	Y_t^{(n)}=\hE\Big[\xi+\int_t^T f(s,Y_s^{(n-1)})ds\ \Big|\  \cF_t\Big],\qq t\in[0,T].
	\eea 
We should point out that  the set-valued random variable $Y_t^{(n)}$ is defined  almost surely for each fixed $t\in [0,T]$. An immediate question is whether $\{Y^{(n)}_t\}_{t\in[0,T]}$ makes sense, as a (jointly) measurable set-valued process, which, as usual, requires justification as we have seen frequently in the set-valued case. The following lemma is important for this 
purpose. 
\begin{lem}\label{lem:progaux}
	Let $X\in \scL^2(\O,\sK(\hR^d))$ and define
$F_t := \hE[X|\cF_t]$,   $t\in[0,T]$.
	Then, $\{F_t\}_{t\in [0,T]}$ has an optional modification that is a uniformly $\hL^2$-bounded martingale.
	\end{lem}

{\it Proof.} Consider the (trivial) set-valued process $G_t\equiv X$,  $t\in[0,T]$, which  is clearly (jointly) measurable and 
 $G_\tau = X$ is integrable for every $\hF$-stopping time $\tau\colon\O\to[0,T]$. By \cite[Theorem 3.7]{optional}, there exists a unique {\it optional projection} $\{{}^{o}G_{t}\}_{t\in [0,T]}$ of process $\{G_t\}_{t\in [0,T]}$, such that
$\hE[G_\tau | \cF_\tau] =  {}^{o}G_{\tau}$, $\hP$-a.s. 
for every $\hF$-stopping time $\tau$. 
In particular, $\{{}^{o}G_{t}\}_{t\in [0,T]}$ is an optional modification of $\{F_t\}_{t\in [0,T]}$. 

It is easy to check that $\{{}^{o}G_{t}\}_{t\in [0,T]}$ is a square-integrable set-valued martingale, and
 by an $\hL^1$-version of Lemma \ref{cond-h2}, it holds that
\beaa
\norm{{}^{o}G_{t}}=\norm{\hE[X|\cF_t]}\leq \hE[\norm{X}|\cF_t], \qq t\in [0,T].
\eeaa
Finally, note that  $X\in \scL^2(\O,\sK(\hR^d))$, applying Doob's $\hL^2$-maximal inequality to the ($\hR$-valued) martingale $M_t:=\hE[\|X\||\cF_t]$, $t\in[0,T]$, we obtain
\beaa
\hE\Big[\sup_{t\in[0,T]}\norm{{}^{o}G_{t}}^2\Big] \leq\hE\Big[\sup_{t\in[0,T]}|M_{t}|^2\Big]\le 4\hE[\norm{X}^2]<+\infty.
\eeaa
That is, $\{{}^{o}G_{t}\}_{t\in [0,T]}$ is uniformly square-integrably bounded.
\qed 

The next proposition establishes the desired measurability for the Picard iteration.

\begin{prop}\label{prop:prog}
	For each $n\in\hN$, $Y^{(n)}$ has a progressively measurable modification. 
\end{prop}

{\it Proof.}
Note that $Y^{(0)}\equiv\{0\}$ is progressively measurable itself. Let $n\in\hN$ and suppose that $Y^{(n-1)}$ has a progressively measurable modification, which we denote by $Y^{(n-1)}$ for ease of notation, and interpret \eqref{picardseq} accordingly.

For each $t\in [0,T]$, using Corollary \ref{condi} and Corollary \ref{dtint-time}, we have
\bea\label{eq:modif}
Y_t^{(n)}= \hE\Big[\xi + \int_0^T f(s,Y_s^{(n-1)})ds\ \Big|\ \cF_t\Big] \ominus \int_0^t f(s,Y_s^{(n-1)})ds.
\eea 
By Remark \ref{propint}(v), $\{\int_0^t f(s,Y_s^{(n-1)})ds\}_{t\in [0,T]}$ has a progressively measurable modification. Moreover, by Lemma \ref{lem:progaux}, $\{ \hE[\xi + \int_0^T f(s,Y_s^{(n-1)})ds |  \cF_t] \}_{t\in [0,T]}$ has an optional, hence progressively measurable, modification. Replacing the original processes with such modifications in \eqref{eq:modif}, and using Lemma \ref{meas}, we have that $Y^{(n)}$ is progressively measurable.
\qed 

In view of Proposition \ref{prop:prog}, we will assume without loss of generality that $Y^{(n)}$ is progressively measurable, in particular, $Y^{(n)}\in \scL^{2}_{\hF}([0,T]\times\O,\sK(\hR^d))$ for each $n\in \hN$.

In order to guarantee the convergence of the sequence $\{Y^{(n)}\}_{n\in\hN}$ constructed in \eqref{picardseq}, we will use a recursive estimate on $\{\hE h^2(Y_t^{(n)},Y_t^{(n-1)})\}_{n\in\hN}$, which is provided by following lemma. We note that unlike the
	vector-valued BSDEs, this lemma is non-trivial because of the lack of standard tools, in particular a set-valued It\^o's formula.

\begin{lem}
	\label{lemma0}
    For each $n\in\hN$, it holds that
	\bea
	\label{DYest1}
	\hE[ h^2(Y_t^{(n)},Y_t^{(n-1)}) ] \leq TK^2\int^T_t \hE [h^2(Y_s^{(n-1)},Y_s^{(n-2)})] ds,\quad t\in[0,T].
	\eea
\end{lem}

	{\it Proof.} By Proposition \ref{norm}, Lemma \ref{cond-h2} and the properties of Hausdorff distance, we get
	\bea
	\hE [h^2(Y_t^{(n)},Y_t^{(n-1)})]
	\negthinspace\negthinspace\negthinspace\negthinspace&=&\negthinspace\negthinspace\negthinspace\negthinspace \hE \sqb{h^2\of{\hE\Big[\xi+\int_t^T f(s,Y_s^{(n-1)})ds \ \Big|\ \cF_t\Big], \hE\Big[\xi+\int_t^T f(s,Y_s^{(n-2)})ds\ \Big|\ \cF_t\Big]}}\nonumber \\
	&\leq &\negthinspace\negthinspace\negthinspace\negthinspace \hE \sqb{\hE\Big[h^2\Big(\xi+\int_t^T f(s,Y_s^{(n-1)})ds,\xi+\int_t^T f(s,Y_s^{(n-2)})ds\Big)\ \Big|\ \cF_t\Big]}\nonumber\\
	&\leq &\negthinspace\negthinspace\negthinspace\negthinspace \hE\sqb{ h^2\of{\int_t^T f(s,Y_s^{(n-1)})ds,\int_t^T f(s,Y_s^{(n-2)})ds}}.\label{holderapp}
	\eea 
	Then, combining \eqref{holderapp} with Proposition \ref{Holder}, Assumption \ref{assump2}(iv), and Proposition \ref{norm}, 
we derive (\ref{DYest1}). 
\qed

We are now ready to establish the well-posedness of the set-valued BSDE \eqref{BSDEc}.

\begin{thm}
	\label{thm:wellposed1}
	Assume Assumption \ref{assump2}. Then, the set-valued BSDE \eqref{BSDEc}
	has a solution $Y\in \scL^{2}_{\hF}([0,T]\times\O,\sK(\hR^d))$. Moreover, the solution is unique up to modifications: if $Y^\prime \in \scL^{2}_{\hF}([0,T]\times\O,\sK(\hR^d))$ is another solution of \eqref{BSDEc}, then $Y_t=Y_t^\prime$ $\hP$-a.s.  $t\in [0,T]$.
\end{thm}

{\it Proof.} 
Recall that $(\scL^{2}_{\hF}([0,T]\times\O,\sK(\hR^d)),d_H)$ is a complete metric space, where the metric is defined by
$d_H(\Phi,\Psi)=(\hE[\int_0^T h^2(\Phi_t,\Psi_t)dt])^\frac12$ for $\Phi,\Psi\in \scL^{2}_{\hF}([0,T]\times\O,\sK(\hR^d))$. We shall argue that the sequence $\{Y^{(n)}\}_{n\in\hN}$ of the Picard iteration is Cauchy in $\scL^{2}_{\hF}([0,T]\times\O,\sK(\hR^d))$. To this end, for fixed $t\in[0,T]$, we note that  $Y_t^{(0)} = \{0\}$. Thus, by repeatedly applying Lemma \ref{cond-h2}, we have
\bea
\label{DY1}
\hE h^2 (Y^{(1)}_t, Y^{(0)}_t)\negthinspace\negthinspace&=&\negthinspace\negthinspace \hE h^2\Big(\hE\Big[\xi + \int^T_t f(s, \{0\}) ds \ \Big|\  \cF_t\Big], \{0\}\Big) \\
&\le&\negthinspace\negthinspace \hE h^2\Big(\xi + \int^T_t f(s, \{0\}) ds, \{0\}\Big) \le 2 \Big[\hE\|\xi\|^2+\hE h^2\Big(\int_t^T f(s, \{0\})ds, \{0\}\Big)\Big]
\nonumber\\
&\le&\negthinspace\negthinspace 2 \Big[\hE\|\xi\|^2+ T\int_0^T\hE \| f(s, \{0\})\|^2ds\Big]=:C.\nonumber
\eea
Note that $C$ is free of the choice of $t$. We claim that, for $n\ge 1$ it holds that
\bea
\label{rec}
\hE h^2 (Y^{(n)}_t, Y^{(n-1)}_t)\leq \frac{C(TK^2)^{n-1}(T-t)^{(n-1)}}{(n-1)!}.
\eea
Indeed, for $n=1$, (\ref{rec}) is just (\ref{DY1}). Now assume that (\ref{rec}) holds for $n-1$, then by Lemma \ref{lemma0}
we have
\bea
\label{DYest5}
\hE h^2 (Y^{(n)}_t, Y^{(n-1)}_t) &\leq& TK^2 \int_t^T \hE\|\D Y^{(n-1)}_s\|^2ds \le TK^2\int^T_t \frac{C(TK^2)^{n-2}(T-s)^{(n-2)}}{(n-2)!}ds\nonumber\\
&=&\frac{C(TK^2)^{n-1}(T-t)^{(n-1)}}{(n-1)!}\le \frac{CK^{2(n-1)}T^{2(n-1)}}{(n-1)!}=: a_n^2.
\eea

Since $\cH_2$ is a metric on $\scL^2_{\cF_t}(\O,\sK(\hR^d))$, the estimate in (\ref{DYest5}) then yields for $m> n\geq 1$:
\bea
\label{YCauchy}
\cH_2(Y^{(n)}_t, Y^{(m)}_t)\leq 
\sum^{m-1}_{k=n}  \cH_2(Y^{(k+1)}_t, Y^{(k)}_t)\le \sum^{m-1}_{k=n}a_{k+1},
\eea
where $a_k= \frac{\sqrt{C}K^{(k-1)}T^{(k-1)}}{\sqrt{(k-1)!}}$, $k\ge 1$, by \eqref{DYest5}. 
Hence,
\bea
\label{YCauchy2}
d^2_H(Y^{(m)},Y^{(n)})&=&\int_0^T \cH_2^2(Y^{(n)}_t, Y^{(m)}_t) dt\le T\Big(\sum^{m-1}_{k=n}a_{k+1}\Big)^2.
\eea
Now note that
\beaa
\frac{a_{k+1}}{a_k} = \frac{\frac{\sqrt{C}K^{k}T^{k}}{\sqrt{k!}}}{\frac{\sqrt{C}K^{(k-1)}T^{(k-1)} }{\sqrt{(k-1)!}}} = \frac{T K}{\sqrt{k}}\to 0, \qq \mbox{\rm
	as $k\to\infty$.}
\eeaa
By ratio test,  $\sum_{k=1}^\infty a_k$ converges.
That is, 
$\{Y^{(n)}\}_{n\in\hN}$ is a Cauchy sequence in $\scL^{2}_{\hF}([0,T]\times\O,\sK(\hR^d))$, thanks to (\ref{YCauchy}); whence converges to some $Y\in\scL^{2}_{\hF}([0,T]\times\O,\sK(\hR^d))$.

Next, we show that the limit process $Y=\{Y_t\}_{t\in[0,T]}$ indeed leads to a solution to the BSDE (\ref{BSDEc}). 
Since $d_H(Y,Y^{(n)})\rightarrow 0$ as $n\rightarrow\infty$, there exists a subsequence $\{Y^{(n_\ell)}\}_{\ell\in\hN}$ such that $h(Y_t(\o),Y^{(n_\ell)}_t(\o))\rightarrow 0$ as $\ell\rightarrow\infty$ for $dt\times d\hP$-a.e. $(t,\o)\in[0,T]\times\O$. By Proposition \ref{norm}(ii) and Corollary \ref{condi}, we have
By Lemma \ref{cond-h2}, Proposition \ref{timeint}, and Assumption \ref{assump2}(iii,iv), we have
\bea
\label{DY2}
&&\hE h^2\Big(\hE\Big[\int^T_t f(s, Y_s) ds\ \Big|\ \cF_t\Big],\hE\Big[\int^T_t f(s, Y^{(n_\ell)}_s) ds\ \Big|\ \cF_t\Big]\Big)\\
&\le& \hE h^2\Big(\int^T_t f(s, Y_s) ds,\int^T_t f(s, Y^{(n_\ell)}_s) ds\Big)\nonumber\\
&\le& (T-t)\hE\Big[\int_t^Th^2(f(s,Y_s),f(s,Y_s^{(n_\ell)}))ds\Big] 
\leq TK^2\int_0^T  \cH_2^2(Y_s,Y_s^{(n_\ell)})ds.\nonumber
\eea
By the construction of the limit $Y$, we have $\int_0^T  \cH_2^2(Y_s,Y_s^{(n_\ell)})ds\rightarrow 0$ as $\ell\rightarrow\infty$.
Now, (\ref{DY2}) shows that 
$$ \sup_{t\in[0,T]}\hE h^2\Big(\hE\Big[\int^T_t f(s, Y_s^{(n_\ell)}) ds\ \Big|\ \cF_t\Big],\hE\Big[\int^T_t f(s, Y_s) ds\ \Big|\ \cF_t\Big]\Big)\to 0, \qq \mbox{as $n\to \infty$.}
$$
It follows that $Y$ satisfies the BSDE
\bea
\label{Y0}  
Y_t=\hE\Big[\xi+\int_t^T f(s, Y_s)ds\ \Big|\ \cF_t\Big], \q t\in [0,T].
\eea
In fact, a similar argument as above (using Assumption \ref{assump2}(iii)) also shows that $Y$ is actually unique, as the solution of  (\ref{BSDEc}) in the space $\scL^{2}_{\hF}([0,T]\times\O,\sK(\hR^d))$ up to modifications. 
This proves the theorem.
\qed

\subsection{Set-Valued BSDEs with Martingale Terms}\label{sec:BSDEmart}

The set-valued BSDE \eqref{BSDEc} considered in Section \ref{sec:condexp} is formulated using set-valued conditional expectations. In this section, we aim to remove the conditional expectation by introducing an additional term to the BSDE, that is, a set-valued martingale. More specifically, we consider a set-valued BSDE of the form
\bea
\label{BSDEmart}
Y_t + M_T = \xi +\int_t^T f(s,Y_s)ds + M_t,\quad t\in[0,T],
\eea
where $\{M_t\}_{t\in[0,T]}$ is a set-valued martingale and $\{Y_t\}_{t\in [0,T]}$, $\xi$, $f\colon [0,T]\times\O\to\sK(\hR^d)$ are as in Section \ref{sec:condexp}. 

\begin{defn}\label{sol:BSDEmart}
	A pair $(Y,M)\in (\scL^{2}_{\hF}([0,T]\times\O,\sK(\hR^d)))^2$ of set-valued process is called an adapted solution to the set-valued BSDE \eqref{BSDEmart} if $M$ is a uniformly square-integrably bounded set-valued martingale with $M_0=Y_0$ and
	\bea
	Y_t + M_T = \xi +\int_t^T f(s,Y_s)ds + M_t,\qq \hP\as, t\in[0,T].\nonumber 
	\eea
\end{defn}

\begin{rem}\label{BSDEmart:Hukuhara}
In light of Hukuhara difference,  the set-valued BSDE \eqref{BSDEmart} is equivalent to
\bea
Y_t = \of{\xi +\int_t^T f(s,Y_s)ds + M_t}\ominus M_T,\quad t\in[0,T].
\eea 
Furthermore, if $M_t\ominus M_T$ exists, the same BSDE is also equivalent to
\bea
Y_t = \xi +\int_t^T f(s,Y_s)ds + (M_t\ominus M_T),\quad t\in[0,T],
\eea 
thanks to Proposition \ref{HD0}. However, the existence of $M_t\ominus M_T$ is a tall order due to the lack of time-additivity of generalized stochastic integral (see Remark \ref{newrem} for more details).
	\end{rem}

We now establish the well-posedness of the set-valued BSDE \eqref{BSDEmart}.

\begin{thm}\label{thm:wellposed2}
	Suppose that Assumption \ref{assump2} holds. Then, the set-valued BSDE \eqref{BSDEmart}
	has a solution $(Y,M)\in (\scL^{2}_{\hF}([0,T]\times\O,\sK(\hR^d)))^2$. Moreover, the solution is unique up to modifications: if $(Y^\prime,M^\prime) \in (\scL^{2}_{\hF}([0,T]\times\O,\sK(\hR^d)))^2$ is another solution of \eqref{BSDEmart}, then $Y_t=Y_t^\prime$ and $M_t=M_t^\prime$ $\hP$-a.s. for every $t\in [0,T]$.
\end{thm}

{\it Proof.} 
By Theorem \ref{thm:wellposed1}, the set-valued BSDE \eqref{BSDEc} in conditional expectation form has a solution $Y\in \scL^{2}_{\hF}([0,T]\times\O,\sK(\hR^d))$. Define a process $M=\{M_t\}_{t\in [0,T]}$ by
\[
M_t:=\hE\Big[\xi+\int_0^T f(s,Y_s)ds\ \Big|\ \cF_t\Big],\quad t\in [0,T].
\]
By Remark \ref{assumpcon},  $\{f(t,Y_t)\}_{t\in [0,T]}\in  \scL^{2}_{\hF}([0,T]\times\O,\sK(\hR^d))$; and by Lemma \ref{lem:progaux}, $M$ is a uniformly square-integrably bounded set-valued martingale. On the other hand, by Corollary \ref{dtint-time},
\bea
\label{add}
\int_0^T f(s,Y_s)ds=\int_0^t f(s,Y_s)ds + \int_t^T f(s,Y_s)ds,\quad t\in[0,T].
\eea 
By the linearity of set-valued conditional expectation and \eqref{add}, we have
\bea
\label{M1}
M_t&\negthinspace\negthinspace=\negthinspace\negthinspace& \hE\Big[\xi + \int^T_0 f(s,Y_s) ds\ \Big|\ \cF_t\Big]= \hE\Big[\xi + \int^t_0 f(s,Y_s) ds + \int^T_t f(s,Y_s) ds\ \Big|\ \cF_t\Big]\\
&=\negthinspace\negthinspace& \hE\Big[\xi + \int^T_t f(s,Y_s) ds\ \Big|\ \cF_t\Big] + \int^t_0 f(s,Y_s) ds= Y_t + \int^t_0 f(s,Y_s) ds. \nonumber
\eea
Using the definitions of $M_T, M_t$, and combining \eqref{add} and \eqref{M1} give
\bea 
Y_t+M_T &=& Y_t + \xi +\int_0^t f(s,Y_s) ds + \int_t^T f(s,Y_s) ds= \xi +\int_t^T f(s,Y_s) ds + M_t.\nonumber \label{M22}
\eea
Finally, $M_0 = \hE[\xi+\int_0^T f(s,Y_s)ds] = Y_0$ by the definitions of $M_0$ and $Y_0$. Hence, the pair $(Y,M)$ is a solution to the set-valued BSDE \eqref{BSDEmart}.

To prove uniqueness, let $(Y^\prime,M^\prime)$ be another solution to \eqref{BSDEmart}. Let $t\in [0,T]$. Hence,
\bea\label{BSDEmartunique}
Y^\prime_t + M^\prime_T = \xi +\int_t^T f(s,Y^\prime_s)ds + M^\prime_t,\qq \hP\as.
\eea 
Since $M^\prime$ is a martingale, taking conditional expectation in \eqref{BSDEmartunique} under $\cF_t$ gives
\bea
Y^\prime_t + M^\prime_t= \hE\Big[\xi+\int_t^T f(s,Y^\prime_s)ds)ds\ \Big| \ \cF_t\Big] + M^\prime_t.
\eea 
Hence, by cancellation law, $Y^\prime$ is a solution to the BSDE \eqref{BSDEc}.
By the uniqueness part of Theorem \ref{thm:wellposed1}, $Y_t=Y_t^\prime$ $\hP$-a.s. Using this, we may rewrite \eqref{BSDEmartunique} as
\bea
Y_t + M^\prime_T = \xi +\int_t^T f(s,Y_s)ds + M^\prime_t,\qq \hP\as.\nonumber 
\eea 
In particular, when $t=0$, we have
\bea\label{BSDEmartunique2}
Y_0 + M^\prime_T = \xi + \int_0^T f(s,Y_s)ds + M^\prime_0.
\eea 
On the other hand, we have $M_0=Y_0=Y_0^\prime=M_0^\prime$. Hence, cancellation law in \eqref{BSDEmartunique2} and the definition of $M_T$ yield
\bea 
M^\prime_T = \xi +\int_0^T f(s,Y_s)ds = M_T.\nonumber 
\eea 
Finally, since $M, M^\prime$ are both martingales with the same terminal value, we obtain $M_t=M^\prime_t$ $\hP$-a.s. for each $t\in [0,T]$. Hence, $(Y,M)$ and $(Y^\prime,M^\prime)$ coincide up to modifications.
\qed

We conclude this section by formulating the precise relationship between the solutions of the two forms of the set-valued BSDE: \eqref{BSDEc} and \eqref{BSDEmart}.

\begin{cor}
	(i) If $Y\in \scL^{2}_{\hF}([0,T]\times\O,\sK(\hR^d))$ is a solution of \eqref{BSDEc}, then there exists a unique $M\in\scL^{2}_{\hF}([0,T]\times\O,\sK(\hR^d))$,
such that $(Y,M)$ is a solution of \eqref{BSDEmart}.

	(ii) If $(Y,M)\in (\scL^{2}_{\hF}([0,T]\times\O,\sK(\hR^d)))^2$ is a solution of \eqref{BSDEmart}, then $Y$ solves \eqref{BSDEc}.
	\end{cor}

{\it Proof.} 
(i) Let $Y$ be a solution of \eqref{BSDEc}. Following the construction in the proof of Theorem \ref{thm:wellposed2}, one can find a set-valued martingale $M$ such that $(Y,M)$ solves \eqref{BSDEmart}. The uniqueness of such $M$ is also a consequence of the uniqueness part of Theorem \ref{thm:wellposed2}.

\ms
(ii) Let $(Y,M)$ be a solution of \eqref{BSDEmart}. For each $t\in[0,T]$, taking conditional expectation with respect to $\cF_t$ in \eqref{BSDEmart} gives
\beaa
Y_t+M_t = \hE\Big[\xi+\int_t^T f(s,Y_s)ds\ \Big| \  \cF_t\Big] + M_t
\eeaa
since $M$ is a set-valued martingale. By cancellation law, it follows that $Y$ solves \eqref{BSDEc}.
\qed

\subsection{Set-Valued BSDEs with Generalized Stochastic Integrals}

We can now combine the martingale representation theorem developed in Section \ref{sec:zero_init} with the set-valued BSDE \eqref{BSDEmart} to study the equation of the form
\bea\label{BSDEint}
Y_t+\int_{0-}^T \cR\circ dB = \xi+\int_t^T f(s,Y_s)ds + \int_{0-}^t \cR\circ dB,\quad t\in [0,T],  
\eea 
where $\cR\subset \hR_0=\hR^d\times \hL^2_{\hF}([0,T]\times\O,\hR^{d\times})$ is a set of martingale representer pairs.

\begin{defn}\label{sol:BSDEint}
	A pair $(Y,\cR)$ with $Y\in \scL^2_{\hF}([0,T]\times\O,\sK(\hR^d))$ and $\cR\in \hR_0$ is called a solution of the set-valued BSDE \eqref{BSDEint} if $Y_0 = \pi_\xi [\cR]$ and
	\bea
	Y_t+\int_{0-}^T \cR\circ dB = \xi+\int_t^T f(s,Y_s)ds + \int_{0-}^t \cR\circ dB,\ \qq \hP\as, t\in[0,T].\nonumber 
	\eea
	\end{defn}

The following theorem provides a well-posedness result for the set-valued BSDE \eqref{BSDEint}.

\begin{thm}\label{thm:wellposed3}
	Under Assumption \ref{assump2}, the set-valued BSDE \eqref{BSDEint} has a solution $(Y,\cR)$. Moreover, the solution is unique in the following sense: if $(Y^\prime,\cR^\prime )$ is another solution of \eqref{BSDEint}, then $Y_t=Y_t^\prime$ and $\int_{0-}^t \cR\circ dB=\int_{0-}^t \cR^\prime\circ dB$ $\hP$-a.s. for every $t\in [0,T]$.
\end{thm}

{\it Proof.} By Theorem \ref{thm:wellposed2}, there exists a solution $(Y,M)$ of the set-valued BSDE \eqref{BSDEmart}. By Definition \ref{sol:BSDEmart}, $M$ is a uniformly square-integrably bounded set-valued martingale. Hence, by Theorem \ref{MRT-initial}, we may write $M_t=\int_{0-}^t \cR^M \circ dB$ $\hP$-a.s. for each $t\in[0,T]$. Hence, $(Y,\cR^M)$ is a solution of \eqref{BSDEint}. The uniqueness claim is an immediate consequence of the uniqueness part of Theorem \ref{thm:wellposed2}.
\qed 

\begin{rem}
	\label{newrem}
	{\rm The indefinite integral $M=\{\int_{0-}^T \cR^M \circ dB\}_{t\in[0,T]}$ in the proof of Theorem \ref{thm:wellposed3} is a uniformly square-integrably bounded set-valued martingale. Following similar arguments as in \cite[Corollary 5.3.2]{Kis2020}, it can be shown that if $\cZ^M$ is decomposable, then, $\cZ^M$ as well as $\cR^M$ are singletons. Hence, in all cases where $\cZ^M$ contains more than one processes, the set $\cZ^M$ is not decomposable. Similar to Corollary \ref{ztimeadd} for the generalized Aumann-It\^o stochastic integral, the indefinite integral considered here does not have time-additivity in general, that is, the Hukuhara difference $M_T\ominus M_t = \int_{0-}^T \cR^M\circ dB\ominus \int_{0-}^t \cR^M\circ dB$ does not exist. In particular, in view of Remark \ref{comparison}, the inclusion $\int_{0-}^t \cR^M\circ dB\subset M_0 +\int_0^t \cZ^M\circ dB$ is generally a strict one.
	\qed}
\end{rem}

\ms

\no{\Large \bf Acknowledgments}

\ms
We would like to thank an anonymous referee who pointed out a serious issue in a previous version of the paper.

  \end{document}